\documentclass[a4paper,11pt,leqno]{amsart}
\usepackage[francais]{babel}
\usepackage[applemac]{inputenc}
\usepackage{epsfig,xypic}
\usepackage{amsthm,amssymb,bbm,a4wide}
\usepackage{hyperref}
\theoremstyle{plain}
\author[I. Marin]{Ivan Marin}
\address{69 rue S\'ebastien Gryphe\\ F-69007 Lyon\\ France}
\date{28 août 2004}
\urladdr{http://www.maths.univ-evry.fr/pages\_perso/~marin}
\title{Monodromie algébrique des groupes d'Artin diédraux}

\newtheorem*{theor}{Th\'eor\`eme}
\newtheorem*{theore}{Th\'eor\`eme {\rm (Enriquez)}}
\newtheorem{lemme}{Lemme}
\newtheorem{definition}{D\'efinition}
\newtheorem*{ccor}{Corollaire}

\newtheorem{prop}{Proposition}

\newcommand{\la}{\lambda}

\def\Ker{\mathrm{Ker}\,}
\def\U{\mathsf{U}}
\def\C{\ensuremath{\mathbbm{C}}}
\def\Q{\mathbbm{Q}}
\def\Z{\mathbbm{Z}}
\def\R{\mathbbm{R}}
\def\N{\mathbbm{N}}

\def\k{\mathbbm{k}}
\def\Ass{\mathbbm{Ass}}

\def\om{\omega}

\def\eps{\epsilon}
\def\g{\mathfrak{g}}

\def\un{\mathbbm{1}}
\def\onto{\twoheadrightarrow}

\def\ii{\mathrm{i}}
\def\kt{\k^{\times}}
\def\bar{\overline}

\def\j{j}
\def\dd{\mathrm{d}}
\def\IXE{\mathcal{X}}
\def\Log{\mathrm{Log}}
\def\det{\mathrm{det}}
\def\cotg{\mathrm{cotg}}
\def\ss{\tilde{\sigma}}
\def\tt{\tilde{\tau}}
\def\sss{\hat{\sigma}}
\def\ttt{\hat{\tau}}
\def\ooo{\hat{O}}
\def\sb{\overline{\sigma}}
\def\tb{\overline{\tau}}
\def\PPhi{\widetilde{\Phi}}
\def\tp{\tau_{\Phi}}
\def\sp{\sigma_{\Phi}}
\def\J{\mathfrak{J}}
\def\Ad{\mathrm{Ad}}
\def\t{\mathrm{t}} 
\def\PP{\mathrm{P}}
\begin{document}

\maketitle

\noindent {\bf Résumé.}
En vue d'étudier la théorie des représentations d'un groupe d'Artin
diédral $B$, nous construisons des morphismes rationnels de
$B$ vers les inversibles de l'algèbre des tresses infinitésimales
associée. Pour ce faire,
nous construisons des
analogues aux associateurs de Drinfeld en rang 2 pour toutes les
symétries diédrales du plan.

\bigskip

\noindent {\bf Abstract.} Towards the study of the representation theory
of any dihedral Artin group $B$, we build rational morphisms
from $B$ to the group of invertible elements of the associated
infinitesimal braids algebra. For this we build analogons of Drinfeld
associators, in rank 2 for arbitrary dihedral symmetries of the plane.

\noindent {\bf MSC 2000.} 20F36, 14F35, 20C08. \\
\noindent {\bf Mots-clés.} Groupes d'Artin, monodromie, associateurs, représentations


\setcounter{section}{-1}
\section{Introduction}

Soit $W$ un groupe de Coxeter fini, et $B$ le groupe de tresses généralisé
ou groupe d'Artin associé. Une présentation de $B$ est déduite de la
présentation de Coxeter de $W$ en n'imposant pas aux générateurs
$\sigma_1,\dots,\sigma_n$ d'être
involutifs. On a donc toujours une projection naturelle $\pi : B \to
W$, dont le noyau $P$ est un sous-groupe caractéristique de $B$ (cf. \cite{COHEN})
appelé groupe des tresses pures
associé à $W$. Cette définition algébrique de $B$ et $P$ est essentiellement
dûe à Tits. Dans le cas où $W$ est de type $A_n$, le groupe d'Artin
associé à $W$ est le groupe de tresses classique sur $n+1$ brins.
Ces groupes naturellement associés aux groupes de Coxeter finis
sont infinis de type fini, sans torsion et linéaires \cite{DIGNE}. Le lien de $B$ avec
$W$ semble suffisamment rigide pour que l'on puisse
espérer esquisser une théorie des représentations de $B$
en lien avec celle de $W$.

Une définition topologique de $B$ et $P$, essentiellement dûe à Brieskorn,
provient de la description de $W$ comme groupe de réflexions
de $\R^n$. Si l'on note $\mathcal{H}$ l'ensemble des hyperplans noyaux des réflexions
de $W$, et $\mathcal{H}_{\C}$ son complexifié dans $\C^n$, $P$ s'identifie
au groupe fondamental de $X = \C^n \setminus \mathcal{H}_{\C}$, et $B$ au groupe
fondamental du quotient de $X$ par l'action naturelle de $W$.

Cet espace $X$ est le complémentaire dans $\C^n$ ou $\mathbbm{P}^n(\C)$ d'un
arrangement central d'hyperplans et un $K(P,1)$. Soit $\mathfrak{g}_X$
son algèbre de Lie (graduée) d'holonomie. Elle admet une présentation simple
sur $\Z$ par générateurs et relations, avec un générateur $\t_H$ (de degré 1)
par élément $H$ de $\mathcal{H}$, et des relations quadratiques. Soit
$\k$ un corps de caractéristique 0, et $\k_h = \k((h))$ le corps
des séries de Laurent sur $\k$.
On note
$\mathfrak{g}_X[\k] = \mathfrak{g}_X \otimes \k$.
L'action
de $W$ sur $X$ induit une action de $W$ sur $\mathfrak{g}_X$, et on a une
1-forme canonique $\Omega$ sur $X$ à valeurs dans $\U \mathfrak{g}_X[\k]$ qui est
$W$-invariante. On note $\U \mathfrak{g}_X(\k)$ la complétion
de $\U \mathfrak{g}_X[\k]$ par rapport à sa graduation naturelle et
$\mathfrak{B}[\k] = \k W \ltimes
\U \mathfrak{g}_X[\k]$, $\mathfrak{B}(\k) = \k W \ltimes
\U \mathfrak{g}_X(\k)$. 
Pour $\k = \C$ les intégrales itérées de Chen
permettent d'associer à chaque choix d'un point-base $z$ de X un morphisme
$\int_z : B \to \mathfrak{B}(\C)$.

L'intérêt de ces morphismes en termes de théorie des représentations
(en dimension finie et en caractéristique 0) est le
suivant. A toute représentation $\check{\rho} : W \to GL_N(\C)$
 de $W$ on peut
associer une variété algébrique
$\mathcal{V}(\check{\rho}) = \{
\rho : \mathfrak{B}[\C] \to M_N(\C) \mid \mathrm{Res}_W \rho =
\check{\rho}\}$
et à tout $\rho \in \mathcal{V}(\check{\rho})$ une représentation
$\rho_h : \mathfrak{B}(\C) \to M_N(\C_h)$ avec $\k_h = \k((h))$,
donc une représentation $\int_z \circ \rho_h : B \to GL_N(\C_h)$ de $B$.
En particulier, chacune des variétés $\mathcal{V}(\check{\rho})$
contient une famille à un paramètre naturelle qui permet d'obtenir les
représentations de l'algèbre d'Iwahori-Hecke générique associée à $W$.
D'autre part, toute représentation irréductible de $\mathfrak{B}[\C]$
donne lieu à une représentation irréductible de $B$ sur $\C_h$.

Enfin, ces représentations sont universellement convergentes en $h$,
et les réprésentations irréductibles obtenues sur $\C_h$ restent
irréductibles après spécialisation en $h \in \C$ pour tout $h$ en dehors d'un ensemble
localement fini de complexes. Comme ces opérations sont compatibles
au coproduit naturel de $\C B$ et de $\mathfrak{B}[\C]$ cela permet
notamment d'introduire des structures naturelle qui reflètent
\og infinitésimalement \fg\ les algèbres d'Iwahori-Hecke,
ou encore de démontrer des résultats d'irréductibilité des
produits tensoriels de représentations de $B$ (cf. \cite{BULL,HECKINF}).

Cependant cette construction de monodromie transcendante, éga\-le\-ment
utilisée dans \cite{BMR} pour la construction d'algèbres de Hecke cyclotomiques,
est entachée de plusieurs défauts. Outre la non-canonicité
du choix d'un point-base dans un espace non contractile, ni même
simplement connexe, elle ne permet
pas de rendre compte de l'apparition fréquente de structures unitaires
sur les représentations de $B$,
et ne garde pas trace du corps de définition de la représentation
infinitésimale $\rho$. En particulier,
il serait utile d'avoir, pour tout corps $\k$ de caractéristique 0, des
morphismes (injectifs) d'algèbres de Hopf
$\PPhi : \k B \to \mathfrak{B}(\k)$ tels que,
posant $\pi(\sigma_i)=s_i$, $H_i = \Ker (s_i - 1)$ et $\t_i = \t_{H_i}$,
$\PPhi$ vérifie la

\medskip

\noindent {\bf Condition fondamentale : } 
Pour tout $i \in [1,n]$, il existe $\la \in \kt$ tel que $\PPhi(\sigma_i)$
est conjugué à $s_i \exp \la \t_i$ par un élément grouplike
de $\U \mathfrak{g}(\k)$.

\medskip

Sous cette hypothèse, les raisonnements
de \cite{ASSOC} s'adaptent im\-mé\-dia\-te\-ment et montrent entre autres
que la correspondance $\widehat{\Phi}$ entre
représentations de $\mathfrak{B}[\k]$ sur $\k$ et représentations
de $B$ sur $\k_h$ naturellement associé à $\PPhi$ est un foncteur
qui préserve toutes les propriétés essentielles de théorie
des représentations. En particulier, on a une notion de
représentations formellement unitaires de $B$ sur $GL_N(\R_h)$
qui permettent d'obtenir, après torsion du corps des coefficients
et spécialisation en $h$ imaginaire pur proche de 0 des
représentations unitaires de $B$. On a également une notion
simple de représentations infinitésimalement unitaires
$\rho : \mathfrak{B}[\R] \to M_N(\R)$, définie par
$\rho(s) \in O_N(\R)$ si $s \in W$ et $\rho(\t_H)$ symétrique, telle que,
si $\rho$ est infinitésimalement unitaire, $\widehat{\Phi}(\rho)$ est
formellement unitaire. L'existence de tels morphismes permettrait notamment
de montrer
qu'il existe sur les représentations des algèbres de Hecke
une structure unitaire directement issue des représentations orthogonales
de $W$. D'autre part, l'exemple du type $A$ montre que les
structures unitaires connues sur les représentations de $B$
ont toutes une interprétation simple au niveau infinitésimal.

On peut construire de tels morphismes pour le type $A$ à l'aide des
associateurs rationnels dont Drinfeld établit l'existence dans \cite{DRIN}.
L'idée
générale de la démonstration de Drinfeld est la suivante. A partir
d'un choix de point-bases tangentiels de $X$ on peut construire
un morphisme $\PPhi : B \to \mathfrak{B}(\C)$, qui vérifie bien 
la propriété voulue pour $\la = \ii \pi$. On introduit alors le
groupe des automorphismes $\Gamma$ de $B$ tels que $\PPhi \circ \Gamma$
vérifie
encore la condition fondamentale. Malheureusement ce groupe discret
est trop petit pour que $\PPhi \circ \Gamma$ puisse être défini sur $\Q$.
Pour pallier à cette difficulté, on montre que $\PPhi$ s'étend en
un morphisme $B(\C) \to  \mathfrak{B}(\C)$,
où $B(\k)$ désigne l'ensemble des $\k$-points d'un schéma en groupe
naturellement relié à la complétion $\k$-pro-unipotente de $P$,
qui contient $B$.
On introduit alors un sous-groupe $G(\k)$ des automorphismes de
$B(\k)$, qui agit librement et transitivement sur l'ensemble de ces
automorphismes. D'un autre côté, les scalaires $\la \in \kt$ agissent
sur $\mathfrak{B}(\k)$ par $\t_H \mapsto \la \t_H$, donc sur l'ensemble
des morphismes considérés. On en déduit qu'à tout $\PPhi$ on peut
associer un sous-groupe à un paramètre de $G(\k)$, donc un élément
non trivial de son algèbre de Lie $\mathcal{G}(\k)$. Si
le groupe $G(\k)$ est choisi suffisamment petit, un tel élément détermine
inversement un unique morphisme $\PPhi$. Il suffit donc de montrer que $\mathcal{G}(\k)$
contient des éléments non triviaux, et cela découle alors du fait
que c'est le cas de $\mathcal{G}(\C)$ à cause de l'existence d'un morphisme
sur $\C$.

C'est cette démonstration que nous adaptons ici pour montrer l'existence
d'un morphisme $\PPhi$ défini sur $\Q$ pour les groupes diédraux,
correspondant au type $I_2(m)$. Pour $m=3$ il s'agit donc d'une
reformulation de la démonstration de Drinfeld, à ceci près que nous
ne prenons pas en compte l'équation dite du pentagone, qui n'a pas de
sens en rang 2. Le groupe $G(\k)$ pour $m=3$ est donc ici plus gros
que le groupe de Grothendieck-Teichmüller considéré dans \cite{DRIN}. Pour
$m$ impair c'est un groupe naturellement associé au groupe d'Artin
$B$, et la démonstration de Drinfeld s'adapte parfaitement. Si
$m$ est pair il faut se restreindre à un sous-groupe de $G(\k)$,
formé des morphismes compatibles avec l'automorphisme non trivial
du diagramme de Coxeter correspondant à $W$
pour que le raisonnement ci-dessus s'applique. Cette distinction
entre les cas $m$ pair et impair recouvre les différences de théorie
des représentations de $W$ entre ces deux cas.

Le résultat principal de cet article peut donc s'énoncer comme suit

\begin{theor} Si $W$ est de type $I_2(m)$, pour tout corps
$\k$ de caractéristique 0 il existe un morphisme $\PPhi :
\k B \to \mathfrak{B}(\k)$ qui vérifie la condition fondamentale.
\end{theor}

Dans \cite{ENRIQUEZ}, B. Enriquez montre l'existence de morphismes en
type $B_n$ qui vérifient la condition fondamentale. De tels
morphismes existent donc au moins en type $A_n$, $B_n$, et $I_2(m)$.
Dans l'appendice 3 nous éclaircissons le lien entre les notations d'Enriquez
et les nôtres.

Une différence notable avec l'approche par intégrales itérées
et point-base est que les morphismes associés, sur un corps $\k \subset
\R$,
ne sont jamais universellement convergents. En effet, si l'un d'entre
eux l'était
toute représentation de $\mathfrak{B}(\R)$ infinitésimalement
unitaire donnerait lieu à des représentations de $B$
unitarisables pour toute spécialisation imaginaire pure de $h$. En
particulier, les représentations de l'algèbre de Hecke $H(q)$ associée
à $W$ seraient unitarisables pour toute valeur de $q$ de module 1.
Comme pour le groupe de tresses habituel ($W$ de type $A_n$),
nous montrons en section \ref{reprW} que ce n'est pas le cas. Inversement, cela suggère que
d'éventuelles propriétés de convergence de ces morphismes sont reliées à des
propriétés de théorie de représentation des groupes d'Artin.

\medskip

\noindent {\bf Conventions générales.} 
Parmi les actions à gauches que nous considérerons, la notation
$g \bullet A$ désignera toujours implicitement la conjugaison par $g$ : 
$g \bullet A = g A g^{-1}$. Comme l'action de $g \in W$ sur $A \in
 \mathcal{A}(\k)$ (définie en \ref{algtressinf})s'identifie à la conjugaison de $A$ par $g$
dans $\mathfrak{B}(\k)$ on utilisera également cette notation
dans ce cadre. On convient d'autre part $g \bullet xy = g \bullet (x y)$
et $g g' \bullet x = (g g') \bullet x$ pour alléger les notations.

Si $(a_i)_{r \leq i \leq s }$ est une famille d'éléments d'un
monoïde, on utilisera les notations
$$
\prod_{i=r}^{s} a_i = a_r a_{r+1} \dots a_s, \ \ 
\prod_{i=s}^{r} a_i = a_s a_{s-1} \dots a_r
$$

Dans ce qui suit, $\k$ désignera un corps quelconque de caractéristique
0, et $\kt$ le groupe
de ses éléments non nuls.
Si $A$ est une $\k$-algèbre, toujours supposée unifère, on désignera par
$A^*$ le groupe de ses éléments inversibles. Si $A$ est une $\k$-algèbre
de Hopf, toujours supposée avec antipode, on désignera par
$A^{\times}$ le groupe de ses éléments grouplike.

\section{G\'en\'eralit\'es}

\subsection{D\'efinitions}

Soit $m \geq 3$, $W$ le groupe di\'edral d'ordre $2m$,
c'est-\`a-dire le groupe de réflexions de type de Coxeter $I_2(m)$.
On note $\theta = \pi/m$ et $c = 2 \cos(\theta)$.
Deux pr\'esentations classiques de $W$ sont donn\'ees par
$
< s, \om \ \mid \ s^2 = \om^m = 1, s \om = \om^{-1} s >$ et
$<s,s' \ \mid \ s^2 = (s')^2 = 1, (s's)^m = 1 >
$,
le passage entre les deux \'etant donn\'e par les
formules $\om = s's$, $s' = \om s$. Les r\'eflexions de $W$ sont
les \'el\'ements de la forme $s \om^r$ pour $0 \leq r \leq m-1$.

On d\'efinit le groupe d'Artin (ou groupe de tresses g\'en\'eralis\'e)
$B$ associ\'e \`a $W$ par la pr\'esentation
$$
< \sigma, \tau \ \mid \ O = \underbrace{\sigma \tau \sigma \dots }_{m\mbox{ facteurs}} =
\underbrace{\tau \sigma \tau \dots }_{m\mbox{ facteurs}} >
$$
On a un morphisme $\pi : B \to W$ d\'efini par $\sigma \mapsto s$,
$\tau \mapsto s'$, dont le noyau $P$ est appel\'e groupe
de tresses pures associ\'e \`a $W$. Le centre de $B$ est cyclique
infini, engendré par $O^2 \in P$ si $m$ est impair
et par $O \not\in P$ si $m$ est pair. La différence essentielle
entre les cas pair et impair au niveau de $W$ tient à ce que, d'une part
le centre de $W$ est cyclique d'ordre 2 engendré par $\om^{\frac{m}{2}}$
si $m$ est pair,
alors qu'il est trivial si $m$ est impair ; d'autre part que $W$ admet deux
classes de conjugaison de réflexions dans le cas pair, une seule dans
le cas impair.

Le groupe $Out(B)$ a été déterminé dans \cite{GILB} (voir également
\cite{PARIS}).
Si $m$ est impair il est cyclique d'ordre 2, engendré
par l'automorphisme $\sigma \mapsto \sigma^{-1}$, $\tau \mapsto \tau^{-1}$
(automorphisme \og image miroir de la tresse \fg), qui est compatible
à la projection $\pi : B \to W$. Si $m$ est pair au contraire
ce groupe est infini, engendré par l'automorphisme image miroir,
un automorphisme d'ordre
infini $\sigma \mapsto \tau^{-1}$, $\tau \mapsto \tau \sigma \tau$,
et l'automorphisme involutif du diagramme de Coxeter $ J : \sigma \mapsto
\tau$, $\tau \mapsto \sigma$. Lorsque $m$ est impair
on peut encore définir $J$, mais il est alors intérieur : il s'agit
de la conjugaison par $O$. Les deux générateurs supplémentaires
de $Out(B)$ dans le cas pair relèvent les automorphismes
de $W$ définis par $s \mapsto \om s$, $\om \mapsto \om$
et $s \mapsto \om s$, $\om \mapsto \om^{-1}$. 
Dans le cas impair
ces automorphismes sont également intérieurs, conjugaisons par
$\om^{\frac{m+1}{2}}$ et $s \om^{\frac{m-1}{2}}$ respectivement.

Tous les
automorphismes de $B$ relèvent donc un automorphisme de $W$,
et en particulier $P$ est un sous-groupe caractéristique de $B$
 pour toute valeur de $m$ (il est stable par tout automorphisme).

\subsection{Algèbre des tresses infinitésimales}
\label{algtressinf}
Soit $m \geq 3$. On pose $\theta = \pi/m$, $\theta_r = r \theta$
et $v_r = (\cos \theta_r, \sin \theta_r)$ pour tout $r \in \Z$. Soit
$\alpha_r \in (\R^2)^*$ la forme linéaire définie par
$\alpha_r(v) = \det(v_r,v)$, $\Delta_r = \Ker \alpha_r \subset \R^2$. Elle
s'étend en une forme linéaire sur $\C^2$ encore notée $\alpha_r$ dont
le noyau est l'hyperplan $H_r \subset \C^2$. Pour toute
$\k$-algèbre (de Lie) $\mathcal{A}$ avec $\k \subset \C$ et
$\t_0,\dots,\t_{m-1} \in \mathcal{A}$, on peut définir la 1-forme
différentielle sur $X = \C^2 \setminus H_0 \cup \dots \cup H_{m-1}$
à valeurs dans $\mathcal{A}$ suivante
$$ \Omega = \sum_{k=0}^{m-1} \t_k \frac{\dd \alpha_k }{\alpha_k}.
$$
On vérifie immédiatement que $\Omega$ est intégrable si et seulement
si
$\t_0 + \dots + \t_{m-1}$ commute à chacun des $\t_r$. L'algèbre
de Lie d'holonomie $\mathfrak{g}$ de $X$, ou algèbre de Lie des tresses infinitésimales
pures associée à $W$, est définie sur $\Q$ par générateurs
$\t_0,\dots,\t_{m-1}$ et relation $T = \t_0 + \dots +\t_{m-1}$ central.
On note $\mathfrak{g}'$ le quotient de $\mathfrak{g}$ par son centre $\Q T$,
et $t_i$ l'image de $\t_i$ dans $\mathfrak{g}'$. Elle
s'identifie à une sous-algèbre de $\mathfrak{g}$ par $t_i \mapsto
\t_i - \frac{T}{m}$. Pour tout corps
$\k$ de caractéristique 0 on note $\mathfrak{g}[\k] =
\mathfrak{g} \otimes \k$, $\mathfrak{g}'[\k] = \mathfrak{g}' \otimes
\k$, $\U\mathfrak{g}[\k] =
\U\mathfrak{g} \otimes \k$, $\U\mathfrak{g}'[\k] = \U\mathfrak{g}' \otimes
\k$. Les algèbres de Lie $\mathfrak{g}[\k]$ et
$\mathfrak{g}'[\k]$ ainsi que leurs algèbres enveloppantes sont
naturellement graduées par $\deg t_i = \deg \t_i = 1$. On note 
$\mathfrak{g}(\k)$, $\mathfrak{g}'(\k)$,
$\mathcal{A}(\k)$ et $\mathcal{A}'(\k)$ leur complétion respectivement
à cette graduation. On notera éventuellement $\t_{H_i} = \t_i$.

Le groupe $W$ agit sur $X$ par la représentation de réflexion,
$s$ agissant par la matrice $\left( \begin{array}{cc} 1 & 0 \\ 0 & -1 \end{array} \right)$
et $\om$ par la rotation d'angle $2 \theta$, et sur $\mathcal{A}(\k)$
par $s.\t_r = \t_{-r}$, $\om . \t_r = \t_{r+2}$, où l'addition est définie
modulo $m$ ($\t_{r+m} = \t_r$). La 1-forme
$\Omega$ à valeurs dans $\mathcal{A}(\k)$ est équivariante pour ces
actions, donc peut être considérée comme une 1-forme sur $X/W$. On
a $P = \pi_1(X)$, $B= \pi_1(X/W)$.

L'action de $W$ sur $\mathcal{A}(\k)$ est héritée de l'action
par conjugaison $g \bullet x = g x g^{-1}$ de $W$ sur l'ensemble de ses réflexions, l'hyperplan
$H_r$ étant stable pour la réflexion $s \om^r$, et $\om \bullet (\om^r s) =
\om^{r+2} s$, $s \bullet( \om^r s) = \om^{-r} s$.
On note $\mathfrak{B}[\k] = \k W \ltimes \U \mathfrak{g}[\k]$
le produit semi-direct d'algèbres de Hopf correspondant, que l'on appelle
algèbre (de Hopf) des tresses infinitésimales associée à $W$, et on
note $\mathfrak{B}(\k) = \k W \ltimes \U \mathfrak{g}(\k)$
sa complétion par rapport à la graduation $\deg(\t_i) = 1$,
$\deg(g) = 0$ si $g \in W$. De façon similaire,
on note $\mathfrak{B}'[\k]= \k W \ltimes \U \mathfrak{g}'[\k]$ le quotient
de $\mathfrak{B}[\k]$ par la relation $T = 0$ et
$\mathfrak{B}'(\k)$ sa complétion.

\subsection{Repr\'esentations de $W$}
\label{reprW}
Un des modèles matriciels de la repr\'esentation de r\'eflexion de $W$
(ou de sa contragr\'ediente suivant la convention choisie)
est donn\'e par
$$
\rho(s) = \left( \begin{array}{cc} -1 & 0 \\ -c & 1 \end{array} \right)
\ \ \ \rho(s') = \left( \begin{array}{cc} 1 & -c \\ 0 & -1 \end{array} \right)
$$
et, de fa\c con g\'en\'erale,
les représentations de $W$ sont définies sur $\k = \Q(\cos \theta) \subset \R$.
De plus, elles peuvent être définies par des matrices
orthogonales sur $\k$. Chacune de ces représentations $\rho : W \to GL_N(\k)$ peut
être étendue à $\mathfrak{B}[\k]$ par $\rho(\t_H) = \rho(s_H)$, où
$s_H \in W$ désigne la réflexion autour de l'hyperplan $H$.
Ces représentations et en particulier l'image par $\rho$ de
$\mathfrak{g}[\k]$ ont été étudiées dans \cite{BULL,HECKINF}.
On en déduit comme dans \cite{ASSOC} que,
une fois construit un morphisme $\PPhi$ vérifiant la condition
fondamentale de l'introduction,
on peut associer à chacune de ces représentations $\rho$ une représentation $\widehat{\Phi}(\rho)
: B \to U_N^{\eps}(\k)$ où $U_N^{\eps}(\k) \subset GL_N(\k_h)$ désigne le groupe
unitaire formel associée à l'automorphisme topologique de
$\k_h$ défini par $h \mapsto -h$. Si un tel morphisme était universellement
convergent, par spécialisation en tout $h$ imaginaire pur on
obtiendrait des représentations unitaires de $B$. Ces
représentations vérifient que l'image de $\sigma$ est semisimple et
a pour spectre
$\{ e^{\la h} , - e^{ -\la h} \}$. On en déduit que, quitte à multiplier
l'image des générateurs par $v = e^{\la h}$, ces représentations
se factorisent par
l'algèbre de Hecke associée au paramètre $v^2$. De plus, si $\rho$ est
irréductible, ces représentations spécialisées en $h$ sont irréductibles
pour toute valeur $h$ en dehors d'un ensemble localement fini de complexes.

L'algèbre de Hecke associée au paramètre $u = v^2$ est le quotient de $\C
B$ par les relations $(\sigma + 1)(\sigma - u) =0$,
$(\tau + 1)(\tau - u) =0$. Posons $d_j = 2 \cos j \frac{\pi}{m}$. Alors
(cf. \cite{GP} lemme 8.1.10, et th. 8.3.1) les représentations irréductibles
de dimension 2 de cette algèbre de Hecke sont toutes de la forme $\rho_j$, pour
$1 \leq j \leq \frac{m-1}{2}$, avec
$$
\rho_j(\sigma) = \left( \begin{array}{cc} -1 & 0 \\ v d_j & v^2
\end{array} \right) \ \ \   
\rho_j(\tau) = \left( \begin{array}{cc} v^2 & v d_j \\ 0 & -1 
\end{array} \right)
$$
Fixons $j$, et soit $a$ (resp. $b$) le projecteur sur l'espace
propre de $\rho_j(\sigma)$ (resp. $\rho_j(\tau)$) correspondant à la
valeur propre $-1$ parallèlement à l'autre espace propre. Si $\rho_j$
est unitarisable en tant que représentation de $B$, cela implique
que $a$ et $b$ sont autoadjoints positifs, donc qu'il en est de même
de $aba$. Mais le spectre de $aba$ est composé de $0$ et de
$d_j^2/(1+2\cos \alpha)$, si $u = \exp \ii \alpha$. Cela n'est donc
possible que si $\cos \alpha > -1/2$. (Cet argument est emprunté
à \cite{WENZL}, prop. 2.9.) En particulier, il
ne peut pas exister de morphisme $\PPhi$ universellement convergent
défini sur $\R$.

\subsection{Complétions de $P$ et de $B$}

Comme $B$ est une extension de $W$ par $P$, $B$ est naturellement isomorphe
\`a $B \ltimes P / \mathcal{Q}$, o\`u $\mathcal{Q}$ est le sous-groupe de $B \ltimes
P$ engendr\'e par les \'el\'ements de la forme $x . x^{-1}$ pour $x
\in P$. Puisque $P$ est une extension centrale d'un groupe libre par $\Z$,
il est résiduellement nilpotent sans torsion, c'est-\`a-dire que
$P$ se plonge dans sa compl\'etion $\k$-prounipotente $P(\k)$
(cf. \cite{BOURBAKI} II \S 6 ex. 4 p.91 et \cite{DELIGNE} p. 175-178). L'action par conjugaison de
$B$ sur $P$ s'\'etend en une action de $B$ sur $P(\k)$ et $B$
se plonge dans $B(\k) = B \ltimes P(\k) / \mathcal{Q}$. On en déduit une surjection
naturelle $\pi : B(\k) \to W$ de noyau $P(\k)$.

Soit $P_r = P/C^r P$ o\`u $C^r P$ est le $r$-i\`eme terme de la suite
centrale descendante de $P$. Le groupe $P_r$ est nilpotent
sans torsion, donc $P_r = \exp L_r$ o\`u $L_r$ est une $\Z$-alg\`ebre de Lie
nilpotente. On a alors $P_r(\k) = \exp L_r \otimes_{\Z} \k$ et
$P(\k)$ est la limite projective des $P_r(\k)$.
C'est un $\k$-sch\'ema en groupe affine ($\k$-groupe) pro-unipotent.
Nous renvoyons à \cite{DG} pour les résultats élémentaires sur ces objets.
Si l'on
d\'esigne par $L(\k)$ la limite projective des $L_r(\k)$, alors
$P(\k) = \exp L(\k)$. On peut alors d\'efinir, pour
tous $X \in P(\k)$, $\la \in \k$, $X^{\la} = \exp \la x$, si $x = \log
X \in L(\k)$.

Si $\Gamma$ est un automorphisme continu de $\k$-groupe de $P(\k)$,
on a alors
$\Gamma(X^{\la}) = \Gamma(\exp(\la x)) = \exp(\la \Gamma(x)) = \Gamma(X)^{\la}$,
o\`u $X = e^x$, $x \in L(\k)$ et $\la \in \k$.

De m\^eme, $B(\k)$ est canoniquement muni d'une
structure de $\k$-groupe. Soit $\Gamma$ un automorphisme de $B(\k)$.
On dit que
$\Gamma$ est \emph{pur} si $\pi \circ \Gamma = \pi$. Il se restreint
alors en un automorphisme continu de $P(\k)$.
Soit $g \in B(\k)$ tel que $\pi(g)$ est d'ordre 2 (par exemple
si $\pi(g)$ est une r\'eflexion), et $\la \in \k$. On d\'efinit
$g^{<\la>} = g (g^2)^{\frac{\la - 1}{2}}$. On a $(g^{<\la>})^{<\mu>}
= g^{<\la \mu>}$. Si $\Gamma$ est pur
alors
$\pi(\Gamma(g)) = \pi(g)$ et
$
\Gamma(g)^{<\la>} = \Gamma(g^{<\la>})$.
On montre alors ais\'ement que
$
\Gamma(g^{<\la>})^{<\mu>} = \Gamma(g)^{<\la \mu>} = \Gamma(g^{<\la \mu>})
$
pour tous $\la,\mu \in \k$.

Supposons donn\'e un morphisme $\Gamma : B \to B(\k)$ tel que $\Gamma(\sigma)
= G^{-1} \sigma^{<\la>} G$, $\Gamma(\tau) = F^{-1} \tau^{<\mu>} F$
avec $F, G \in P(\k)$ et $\la,\mu \in \k$. Il s'\'etend
suivant les inclusions naturelles $L_r \hookrightarrow L_r \otimes_{\Z} \k$
en un automorphisme de $\k$-groupe de $B(\k)$ encore not\'e $\Gamma$
tel que $\pi \circ \Gamma = \pi$.

Le centre de $P$ est engendr\'e par
$Z = (\sigma \tau)^m = (\tau \sigma)^m$, et $P$
lui-m\^eme est engendr\'e par $u_0,\dots,u_{m-1}$ o\`u l'on a
d\'efini
$$
u_0 = \sigma^2, u_1 = \tau^2, u_{2r} = (\tau \sigma)^r \bullet \sigma^2,
u_{2r+1} = (\tau \sigma)^r \bullet \tau^2.
$$
On a $u_{m-1} \dots u_1 u_0 = (\sigma \tau)^m = Z$ et, \`a partir
de la description g\'eo\-m\'e\-tri\-que de $P$ comme groupe fondamental
on montre qu'une pré\-sen\-ta\-tion de $P$ est
$$
<u_0,\dots,u_{m-1},Z \ \mid \ u_{m-1} \dots u_1 u_0 = Z, (Z,u_i) = 1 >
$$
o\`u $(x,y) = xyx^{-1} y^{-1}$. On en d\'eduit que $L(\Q)$
est topologiquement engendr\'ee par des g\'en\'erateurs
$v_0,\dots,v_{m-1},z$ avec $u_i = \exp v_i$, $Z = \exp z$ et
relations $v_0 + \dots + v_{m-1} = z$, $[z,v_i] = 0$. D'autre
part, on a un morphisme $P \to <Z>$ d\'efini par $u_i \mapsto Z$,
$Z \mapsto Z^m$ qui s'\'etend en un morphisme $P(\Q) \to \exp \Q z$
dont la diff\'erentielle $L(\Q) \to \Q z$ est $v_i \mapsto z$, $z \mapsto mz$.
Ceci permet de l'\'etendre en un morphisme $P(\k) \to \exp \k z$. Ce
morphisme est scind\'e par $Z \mapsto Z^{\frac{1}{m}} = \exp (\frac{z}{m})$,
ce qui permet de d\'ecomposer $P(\k)$ en un produit direct $P'(\k)
\times \exp \k z$.

\section{L'associateur transcendant}
\setcounter{equation}{0}

\subsection{R\'eduction \`a une variable}

On rappelle que $\theta = \pi/m$, $\theta_r = r\theta$,
$v_r = (\cos \theta_r,\sin \theta_r)$,
$\alpha_r(v) = \det(v_r,v)$. On a $\Ker \alpha_r = \R v_r = \Ker \alpha_{r+m}$
o\`u l'addition est modulo $2m$. Pour $0 \leq r \leq 2m-1$, on note
$U_r = \R_+^* v_r + \R_+^* v_{r+1}$ le demi-c\^one ouvert engendr\'e par
$v_r$ et $v_{r+1}$, $\Delta_r = \R v_r$, $\Delta = \bigcup_{r=0}^{m-1}
\Delta_r$.

En coordonn\'ees polaires, tout $v \in \R^2 \setminus \{ 0 \}$
s'\'ecrit uniquement sous la forme $v = (\rho \cos u, \rho \sin u)$
avec $\rho > 0$, $u \in [0,2\pi[$. On a alors
$\alpha_r(v) = \rho \sin(u - \theta_r)$ et $U_r$ s'identifie
(analytiquement) \`a $\R_+^* \times V_r$, o\`u $V_r = ]\theta_r,
\theta_{r+1}[$. On consid\`ere le syst\`eme diff\'erentiel
sur $\R^2 \setminus \Delta$ \`a valeurs dans $\mathcal{A}(\R)$

\begin{eqnarray} \label{syst1}
\dd F = \Omega F, \ \ \ \Omega = \sum_{r=0}^{m-1} \t_r \frac{\dd \alpha_r}{\alpha_r}.
\end{eqnarray}

Comme $\dd \alpha_r = \sin(u-\theta_r) \dd \rho + \rho \cos(u-\theta_r)
\dd u$, on a
$$
\Omega = \sum_{r=0}^{m-1} \t_r \left(
\frac{\dd \rho}{\rho} + \cotg(u-\theta_r) \dd u \right)
$$
o\`u $\cotg z = \cos z / \sin z$. On en d\'eduit que (\ref{syst1})
est \'equivalent \`a 

\begin{eqnarray} \label{syst2}
\frac{\partial F}{\partial \rho} = \frac{T}{\rho} F \ \ \ 
\frac{\partial F}{\partial u} = \left( \sum_{r=0}^{m-1} \t_r
\cotg(u-\theta_r) \right) F
\end{eqnarray}
avec $T = \sum_{r=0}^{m-1} \t_r$. Il est donc \'equivalent de consid\'erer
le syst\`eme \`a valeurs
dans $\mathcal{A}'(\R)$ donn\'e par $\partial F/ \partial \rho = 0$ et

\begin{eqnarray} \label{syst3}
\frac{\partial F}{\partial u} = \left( \sum_{r=0}^{m-1} t_r
\cotg(u-\theta_r) \right) F
\end{eqnarray}
Ainsi $F$ est une fonction de
$u$, et du fait que $\cotg(u) - 1/u$ est analytique en 0 on d\'eduit
(cf. appendice 1, lemmes \ref{unicitepositive} et
\ref{unicitenegative}) qu'il existe des solutions uniques de
(\ref{syst3}) sur $V_r$ telles que
$$
F_{r,+} \sim (u-\theta_r)^{t_r}, \ u \to \theta_r^+, \ \ F_{r+1,-} \sim
(\theta_{r+1}-u)^{t_{r+1}}, \ u \to \theta_{r+1}^-,
$$
o\`u l'on a pos\'e $t_{r+m} = t_r$,
qui sont \`a valeurs grouplike (cf. appendice 1, lemme \ref{grouplike}).

On consid\`ere l'\'equation diff\'erentielle (\ref{syst3}) pour
$u$ une variable complexe, c'est-\`a-dire d'abord
$u \in \mathcal{B} = \left( [0,2\pi[  + \ii \R \right) \setminus \Theta$, o\`u
$\Theta = \frac{\pi}{m} \Z$. On note $D_r = \theta_r + \ii \R_-$,
$\mathcal{B}' = \mathcal{B} \setminus \bigcup_{r=0}^{m-1} D_r$.
Comme $\mathcal{B}'$ est simplement connexe, chacune des fonctions
$F_{r,\pm}$ s'\'etend en une fonction analytique sur $\mathcal{B}'$.
On en déduit (cf. appendice 1, lemme \ref{passage}) que $F_{r,+} = F_{r,-}
\exp(\ii \pi t_r)$.

\subsection{\'Equation du demi-tour}

L'\'equation (\ref{syst3}) s'\'ecrit $\dd F = \Omega' F$,
o\`u $\Omega' = \sum_{r=0}^{m-1} t_r \cotg(u -\theta_r) \dd u$
est une 1-forme sur $C = (\C \setminus \Theta)/ \pi \Z$. Notons
$\widehat{C} = C \cup \{ \ii \infty, -\ii \infty \}$ le
compactifi\'e de $C$ en $\ii \infty$ et $- \ii \infty$,
$\widehat{\C} = \C \cup \{ \infty \}$ la sph\`ere de Riemann,
$\zeta = \exp(2 \ii \theta)$ et $\mu_m = \{ \zeta^r \ \mid \ r \in \Z \}$.
On a un isomorphisme $ \widehat{C} \setminus \Theta \to
\widehat{C} \setminus \mu_m$ donn\'e par l'application
$\varphi : u \mapsto \exp(2 \ii u)$, avec $\varphi(\ii \infty) = 0$,
$\varphi(- \ii \infty) = \infty$. On a $\cotg(u) = (\varphi(u)
+ 1)/(\varphi(u) - 1)$, donc
$$
\varphi_* \Omega' = \sum_{r=0}^{n-1} t_r
\frac{z \zeta^{-r} + 1}{z \zeta^{-r} - 1} \frac{\dd z}{z}
= \sum_{r=0}^{m-1} t_r \frac{z+\zeta^r}{z - \zeta^r} \frac{\dd z}{z}
= \sum_{r=0}^{m-1} t_r \frac{2\zeta^r}{z-\zeta^r} \frac{\dd z}{z}
$$
est une 1-forme sur $\C \setminus \left( \mu_m \cup \{ 0 \} \right)$.
Comme $\sum_{r=0}^{m-1} t_r = 0$, on montre ais\'ement que
$\varphi_* \Omega'$ s'\'etend en une 1-forme ferm\'ee sur $\widehat{\C}\setminus
\mu_m$, donc que $\Omega'$ s'\'etend en une 1-forme ferm\'ee sur
$\widehat{C}
\setminus \Theta$. En particulier, le chemin emprunt\'e pour
prolonger analytiquement $F_{0,+}$ au voisinage de $\pi^+$ devient un
lacet homotopiquement trivial de $\widehat{C} \setminus \Theta$,
d'o\`u $F_{0,+}(u) = F_{0,+}(u+ \pi)$,
donc $F_{0,+}(u) \sim (u- \pi)^{t_0}$ quand $u \to 0$
et $F_{0,+} = F_{m,+}$.

\begin{figure}
\begin{center}
\includegraphics{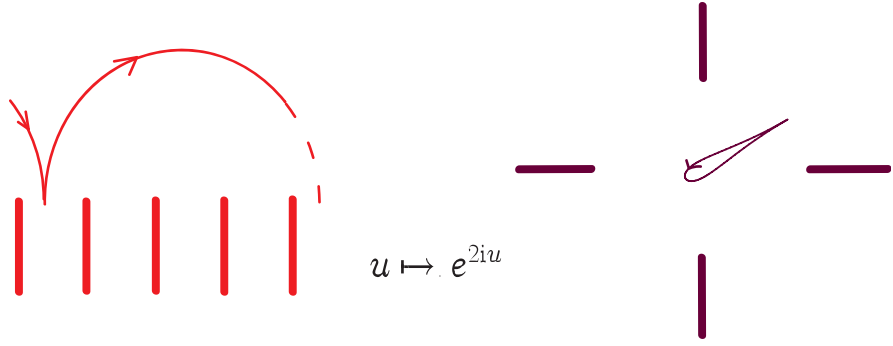}
\end{center}
\caption{Trivialit\'e homotopique du demi-tour}
\end{figure}

\subsection{Action du groupe di\'edral}

On con\-si\-dère \`a nou\-veau l'\'e\-qua\-tion (\ref{syst3}), cette fois sur
$Y = (\R \setminus \Theta)/\pi \Z$. On a une action naturelle
du groupe dié\-dral
$W$ sur $\R^2$, o\`u $s$ agit par la matrice
$\left( \begin{array}{cc} 1 & 0 \\ 0 & -1 \end{array} \right)$
et $\om$ par la rotation d'angle $2 \theta$. Si $v = (\rho \cos u,
\rho \sin u) \in \R^2 \setminus \{ 0 \}$, on a
$$
\begin{array}{lcl}
s . v & = & (\rho \cos(-u), \rho \sin(-u) )\\
\om . v & = & (\rho \cos(u+2\theta), \rho \sin(u+ 2 \theta) )
\end{array}
$$
donc on en d\'eduit une action de $W$ sur $Y$, par $s . u
= -u$, $\om . u = u + 2 \theta$. De l'action de $W$ sur $\mathcal{A}'(\R)$
on déduit une action de $W$ sur le faisceau $\mathcal{F}$
des solutions de (\ref{syst3}) sur $Y$. En effet, soit $\mathcal{G}$ le
faisceau des fonctions analytiques sur $Y$. Si $U$ est un
ouvert de $Y$, $f \in \mathcal{G}(U)$ et $g \in W$ on note
$g . f \in \mathcal{G}(g . U) : u \mapsto g .  f(g^{-1} . u)$. \'Etendant
cette action aux formes diff\'erentielles, on remarque que
$g . \Omega' = \Omega'$ pour tout $g \in W$, donc
si $f \in \mathcal{F}(U)$ on a $g . f \in \mathcal{F}(g .U)$.
Par comparaison des comportements asymptotiques on en d\'eduit
$\om \bullet F_{r,\pm} = F_{r+2,\pm}$
et
$s \bullet F_{r,\pm} = F_{2m-r,\mp} = F_{-r,\mp}$.

Par le th\'eor\`eme de Cauchy sur $V_0$ et $V_1$, il existe
$\Phi_0, \Phi_1 \in \mathcal{A}'(\R)$ grouplike tels que
$F_{1,-} = F_{0,+} \Phi_0$, $F_{2,-} = F_{1,+} \Phi_1$.
On en d\'eduit
\begin{eqnarray}
F_{1,+} = F_{1,-} \exp(\ii \pi t_1) = F_{0,+} \Phi_0 \exp(\ii \pi t_1) 
\label{eqa} \\
F_{2,+} = F_{2,-} \exp(\ii \pi t_2) = F_{1,+} \Phi_1 \exp(\ii \pi t_2)
\label{eqb}
\end{eqnarray}

Appliquant $\om^r$ \`a ces deux \'egalit\'es, on a
\begin{eqnarray} \label{eqc}
\left\lbrace \begin{array}{l}
F_{2r+1,+} =  F_{2r,+} (\om^r \bullet \Phi_0) \exp(\ii \pi t_{2r+1}) \\
F_{2r+2,+} =  F_{2r+1,+} (\om^r \bullet \Phi_1) \exp(\ii \pi t_{2r+2})
\end{array}
\right.
\end{eqnarray}
Appliquant $\om s$ \`a l'\'equation $F_{1,-} = F_{0,+} \Phi_0$ il vient
$F_{1,+} = F_{2,-} (\om s \bullet \Phi_0)$. Comme $F_{2,-} = F_{1,+}
\Phi_1$ on obtient
\begin{eqnarray}
\Phi_1^{-1} = \om s \bullet \Phi_0 \label{phi0phi1}
\end{eqnarray}
On note $\PP_0 = \Phi_0 \exp(\ii \pi t_1)$, et
$$
\xi = \Phi_0 \exp(\ii \pi t_1) \Phi_1 \exp(\ii \pi t_2 )
= \Phi_0 \exp(\ii \pi t_1) (\om s \bullet \Phi_0^{-1}) \exp(\ii \pi t_2).
$$
Comme $F_{2,+} = F_{0,+} \xi$, on déduit de l'action de $\om^r$
que $F_{2r+2,+} = F_{2r,+} (\om^r  \bullet \xi)$ d'o\`u, pour tout $r \geq 0$,
\begin{eqnarray} \label{eqcycle}
F_{2r+2,+} = F_{0,+} \prod_{j=0}^r (\om^j \bullet \xi)
\end{eqnarray}

\subsection{Cas $m$ impair} On suppose que $m \geq 3$ est impair.
\begin{lemme} \label{debeqtr} $\Phi_0$ et $\Phi_1$ v\'erifient les \'equations
$\Phi_0 = \om^{\frac{m-1}{2}} \bullet \Phi_1$, $\Phi_1 = \om^{\frac{m+1}{2}}
\bullet \Phi_0$, $\Phi_0^{-1} = s \om^{\frac{m-1}{2}} \bullet \Phi_0$,
$\Phi^{-1} = s \om^{\frac{m-3}{2}} \bullet \Phi_1$.
\end{lemme}
\begin{proof}
Ces quatre \'equations se d\'eduisent toutes de l'une d'entre
elles par action de $W$ et la relation (\ref{phi0phi1}). Il suffit
donc de d\'emontrer $\Phi_0 = \om^{\frac{m-1}{2}} \bullet \Phi_1$.
Appliquant $\om^{\frac{m-1}{2}}$ \`a (\ref{eqb}) et en utilisant $F_{m+1,+}
= F_{1,+}$, $F_{m,+} = F_{0,+}$, il vient 
$F_{1,+} = F_{0,+} (\om^{\frac{m-1}{2}} \bullet \Phi_1) \exp(\ii \pi t_1)$.
De (\ref{eqa}) on d\'eduit alors $\Phi_0 = \om^{\frac{m-1}{2}} \bullet
\Phi_1$.
\end{proof}
On a alors
$$
\xi = \Phi_0 \exp(\ii \pi t_1) ( \om s \bullet \Phi_0^{-1} ) \exp(\ii \pi t_2)
= \PP_0 ( \om^{\frac{m+1}{2}} \bullet \PP_0 )
$$
\begin{lemme}
$\Phi_0$ v\'erifie l'\'equation
$$
\left( \prod_{r=0}^{\frac{m-3}{2}} \om^r \bullet \xi \right) (\om^{\frac{m-1}{2}}
\bullet \PP_0) = 1
$$
c'est-\`a-dire
$$
\left( \prod_{r=0}^{\frac{m-3}{2}}  (\om^r \bullet \PP_0)
(\om^{\frac{m+1}{2}+r} \bullet \PP_0) \right) (\om^{\frac{m-1}{2}}
\bullet \PP_0) = 1
$$
\end{lemme}
\begin{proof} (\ref{eqc}) se r\'e\'ecrit $F_{2r+1,+} = F_{2r,+} (\om^r \bullet
\PP_0)$ et en particulier $F_{m,+} = F_{m-1,+}(\om^{\frac{m-1}{2}} \bullet
\PP_0)$. L'\'enonc\'e se d\'eduit alors de $F_{m,+} = F_{0,+}$ et de la
relation (\ref{eqcycle}) pour $r = (m-3)/2$.
\end{proof}

\subsection{Cas $m$ pair} On suppose que $m \geq 4$ est pair.

On d\'eduit comme dans le cas $m$ impair de la relation (\ref{eqcycle}),
avec cette fois $r = \frac{m}{2} - 1$, et de l'invariance
par demi-tour $F_{m,+} = F_{0,+}$ la relation suivante
\begin{lemme} \label{fineqtr} On a
$$
\prod_{r=0}^{\frac{m}{2} -1} (\om^r \bullet \xi) = 1
$$
\end{lemme}
Nous verrons en section \ref{section6}, lemme \ref{lemmtranspair}, une autre relation v\'erifi\'ee par
$\Phi_0$, qui fait intervenir un automorphisme \emph{ext\'erieur}
de $\mathfrak{B}(\k)$. Cet élément $\Phi_0$ vérifie
d'autres relations algébriques (par exemple la relation du
pentagone de Drinfeld si $m=3$), étroitement liées aux propriétés
de dépendance algébrique sur $\Q$ de ses coefficients. En particulier,
$\Phi_0$ et ses coefficients semblent reliés aux séries formelles
et aux nombres zêtas multiples
qui apparaissent dans \cite{ENRIQUEZ}. Cependant, seules
les relations établies ici sont utiles pour définir les
morphismes rationnels qui nous intéressent dans ce travail.

\section{Associateurs di\'edraux}
\setcounter{equation}{0}

\subsection{D\'efinitions}

Soit $\k$ un corps de caractéristique 0. L'algèbre $\mathcal{A}'(\k)$
est une algèbre de Hopf complète au sens de Quillen (cf. \cite{QUILLEN}
appendice A2). On note $\widehat{\otimes}$ le produit tensoriel
complété et $\Delta : \mathcal{A}'(\k) \to \mathcal{A}'(\k)\widehat{\otimes}
\mathcal{A}'(\k)$
son coproduit. Si $x,y \in \mathcal{A}'(\k)$ on identifie
$x \otimes y$ à son image dans  $\mathcal{A}'(\k)\widehat{\otimes}
\mathcal{A}'(\k)$.

\begin{definition} Soit $m \geq 3$ impair. Pour tout $\la \in \k$,
on note $\Ass^{(m)}_{\la}(\k)$ l'ensemble des $\Phi \in \mathcal{A}'(\k)$ tels que
\begin{eqnarray}
\Delta(\Phi) = \Phi \otimes \Phi \label{impeq1} \\
\Phi^{-1} = s \om^{\frac{m-1}{2}} \bullet \Phi \label{impeq2} \\
\left[ 
\prod_{r=0}^{\frac{m-3}{2}} (\om^r \bullet \PP) (\om^{\frac{m+1}{2} + r}
\bullet \PP) \right] (\om^{\frac{m-1}{2}} \bullet \PP) = 1
\label{impeq3}
\end{eqnarray}
avec $\PP = \Phi e^{\la t_1}$.
\end{definition}

\begin{definition} Soit $m \geq 4$ pair. Pour tout $\la \in \k$,
on note $\Ass^{(m)}_{\la}(\k)$ l'ensemble des $\Phi \in \mathcal{A}'(\k)$ tels
que
\begin{eqnarray}
\Delta(\Phi) = \Phi \otimes \Phi \label{paireq1} \\
\prod_{r=0}^{\frac{m}{2}-1} (\om^r \bullet \xi)  = 1
\label{paireq3}
\end{eqnarray}
avec $\xi = \Phi e^{\la t_1} (\om s \bullet \Phi^{-1}) e^{\la t_2}$.
\end{definition}

Que $m$ soit pair ou impair, l'action naturelle $W$-invariante
de $\kt$ sur $\mathcal{A}(\k)$ telle que
$\mu \bullet \t_r = \mu \t_r$ laisse stable $\mathcal{A}'(\k)$
et envoie surjectivement $\Ass^{(m)}_{\la}(\k)$ sur
$\Ass^{(m)}_{\la \mu}(\k)$. D'autre part, les \'el\'ements grouplike $\Phi_0$
construits de fa\c con transcendante appartiennent \`a $\Ass^{(m)}_{\ii \pi}(\C)$,
d'apr\`es les lemmes \ref{debeqtr} \`a \ref{fineqtr}, donc $\Ass^{(m)}_1(\C) \neq 0$ pour tout $m \geq 3$. Remarquons que pour
$m = 3$ on retrouve les \'equations des associateurs
de Drinfeld, l'\'equation du demi-tour (\ref{impeq3}) \'etant
celle dite de l'hexagone (l'\'equation du pentagone n'a pas de
sens en rang 2). 

Au premier ordre, on peut d\'ecrire $\Phi \in \Ass^{(m)}_{\la}(\k)$
par $\Phi \equiv 1 + \sum_{r=0}^{m-1} a_r t_r$ avec $a_r \in \k$ et on
peut supposer $\sum_{r=0}^{m-1} a_r = 0$. Si $m$ est impair,
l'\'equation (\ref{impeq2}) impose $a_{m-r+1} = -a_r$ pour tout
$r$, où on a convenu $a_{m+r} = a_r$. En particulier, $a_1 = -a_0$.
Dans tous les cas, le lemme suivant permet d'\'eliminer les coefficients
$a_0$ et $a_1$, en prenant $f(X) = e^{a_0 X}$. Il refl\`ete,
dans le cas de l'associateur transcendant $\Phi_0$, le fait que
l'on aurait pu modifier les conditions asymptotiques sur $F_{0,+}$
(resp. $F_{1,-}$) par une fonction grouplike de $t_0$ (resp. de $t_1$).
C'est un analogue du changement de jauge ({\it twist}) de Drinfeld.

\begin{lemme} Soit $f(X) = \exp(\alpha X)$ avec $\alpha \in \k$
et $\Phi \in \Ass^{(m)}_{\la}(\k)$. Alors $\Phi' =
f(t_0)^{-1} \Phi f(t_1) \in \Ass^{(m)}_{\la}(\k)$. De plus, si $m$ est pair,
$\Phi f(t_1) \in \Ass^{(m)}_{\la}(\k)$.
\end{lemme}
\begin{proof}
On a bien $\Delta(\Phi') = \Phi' \otimes \Phi'$. Si $m$ est impair,
$$
s \om^{\frac{m-1}{2}} \bullet \Phi' = f(t_1)^{-1}
(s \om^{\frac{m-1}{2}} \bullet \Phi) f(t_0) =
f(t_1)^{-1} \Phi^{-1} f(t_0) = (\Phi')^{-1}
$$
d'o\`u (\ref{impeq2}) et, si l'on pose $\PP' = \Phi' e^{\la t_1}$,
$$
\PP' = f(t_0)^{-1} \Phi f(t_1) e^{\la t_1} = f(t_0)^{-1} \PP f(t_1).
$$ 
Que $m$ soit pair ou impair, si $\xi = \Phi e^{\la t_1} (\om s \bullet
\Phi^{-1}) e^{\la t_2}$ et $\xi' = \Phi' e^{\la t_1} (\om s
\bullet (\Phi')^{-1}) e^{\la t_2}$, on a
$$
\xi' = f(t_0)^{-1} \Phi f(t_1) e^{\la t_1} f(t_1)^{-1} (\om s \bullet
\Phi^{-1}) f(t_2) e^{\la t_2}
$$
soit $\xi' = f(t_0)^{-1} \xi f(t_2)$ et $\om^r \bullet \xi' =
f(t_{2r})^{-1} (\om^r \bullet \xi) f(t_{2r+2})$. On d\'eduit alors
de (\ref{impeq3}) et (\ref{paireq3}) que $\Phi' \in \Ass^{(m)}_{\la}(\k)$,
soit la premi\`ere partie de l'\'enonc\'e.
Si $m$ est pair et $\Phi' = \Phi f(t_1)$, $\Phi'$ est
grouplike, $\om s \bullet \Phi' = (\om s \bullet \Phi) f(t_1)$
et $\xi' = \xi$ d'o\`u $\Phi' \in \Ass^{(m)}_{\la}(\k)$. 
\end{proof}
 
La proposition suivante peut se d\'emontrer par un calcul direct,
assez long. Elle peut \'egalement se d\'eduire de
l'automorphisme \og image miroir \fg\ de $B$. C'est ce que nous ferons
à la fin de cette partie.
\begin{prop} \label{assplusmoins}
Pour tout $\la \in \k$, $\Ass^{(m)}_{\la}(\k) = \Ass^{(m)}_{-\la}(\k)$.
\end{prop}

\subsection{Relations entre $\tilde{\sigma}$ et $\tilde{\tau}$}
On fixe $\Phi \in \Ass^{(m)}_{\la}(\k)$, et on pose
$$
\ss = s e^{\la t_0}, \ \ \ \ \tt = \Phi \om s e^{\la t_1} \Phi^{-1}.
$$
On a
$$
\begin{array}{lcl}
\tt \ss & = & \Phi \om s e^{\la t_1} \Phi^{-1} s e^{\la t_0} \\
& = & \om s (\om s \bullet \Phi) e^{\la t_1} \Phi^{-1} s e^{\la t_0} \\
& = & \om (\om^{-1} \bullet \Phi) e^{\la t_{-1}} (s \bullet \Phi)^{-1}
e^{\la t_0} \\
& = & \om ( \om^{-1} \bullet \xi) 
\end{array}
$$
d'o\`u, pour tout $r \geq 1$,
\begin{eqnarray} \label{ttss}
(\tt \ss)^r = \om^r \prod_{j=r}^1 \om^{m-j} \bullet \xi = 
\om^r \prod_{q=m-r}^{m-1} (\om^q \bullet \xi)
\end{eqnarray}

\begin{lemme} \label{relationssttimp} Si $m$ est impair, $(\tt \ss)^{\frac{m-1}{2}} = 
\om^{\frac{m-1}{2}} \PP^{-1}$ et
$$
\ss (\tt \ss)^{\frac{m-1}{2}} = s \om^{\frac{m-1}{2}} \Phi^{-1}
= (\tt \ss)^{\frac{m-1}{2}} \tt
$$
\end{lemme}
\begin{proof} Calcul direct \`a partir de (\ref{ttss}).
\end{proof} 

\begin{lemme} \label{relationssttpair} Si $m$ est pair, $(\tt \ss)^{\frac{m}{2}}
= (\ss \tt)^{\frac{m}{2}} = \om^{\frac{m}{2}}$.
\end{lemme}
\begin{proof}
On d\'eduit par calcul direct de (\ref{ttss}) que $(\tt \ss)^{\frac{m}{2}} =
\om^{\frac{m}{2}}$.
D'autre part, $\ss^{-1} = s e^{-\la t_0}$, $\tt^{-1} = 
\Phi \om s e^{-\la t_1} \Phi^{-1}$. Comme $\Phi \in \Ass^{(m)}_{\la}(\k)
= \Ass^{(m)}_{-\la}(\k)$, on en d\'eduit que
$ (\tt^{-1} \ss^{-1})^{\frac{m}{2}} = \om^{\frac{m}{2}}$,
d'o\`u $(\ss \tt)^{\frac{m}{2}} = \om^{-\frac{m}{2}} =
\om^{\frac{m}{2}}$.
\end{proof}

\subsection{Aspects particuliers du cas $m$ pair}

Au niveau des groupes di\'edraux, les principales diff\'erences
entre les cas $m$ pair et impair proviennent du fait que, si $m$
est pair, $W$ admet deux classes de conjugaisons de r\'eflexions,
et non plus une. Ceci a pour cons\'equence que, sous l'action de $W$,
les g\'en\'erateurs $\t_r$ de $\mathfrak{g}(\k)$ se
r\'epartissent en deux orbites, suivant que $r$ est pair ou impair. En
particulier, on peut \'etendre l'action de $\kt$ sur $\mathcal{A}
(\k)$ en une action de $\kt \times \kt $,
par $(\alpha ,\beta) . \t_r = \alpha \t_r$ si $r$ est pair,
$(\alpha,\beta) . \t_r = \beta \t_r$ si $r$ est impair. On peut alors
d\'efinir, pour $\la, \mu \in \k$, l'ensemble $\Ass^{(m)}_{\la,\mu}(\k)$
des $\Phi \in \mathcal{A}'(\k)$ v\'erifiant les \'equations
(\ref{paireq1}) et (\ref{paireq3}), avec
cette fois $\xi = \Phi e^{\la t_1} (\om s \bullet \Phi^{-1})
e^{\mu t_2}$. On voit facilement que $(\alpha,\beta) .
\Ass^{(m)}_{\la,\mu}(\k) = \Ass^{(m)}_{\alpha \la,\beta\mu}(\k)$, 
$\Ass^{(m)}_{\la,\la}(\k) = \Ass^{(m)}_{\la}(\k)$.
Enfin,
si $\Phi \in \Ass^{(m)}_{\la,\mu}(\k)$ et $\ss = s e^{\la t_0}$,
$\tt = \Phi \om s e^{\mu t_1} \Phi^{-1}$, on montre de m\^eme
$(\ss \tt)^{\frac{m}{2}} = (\tt \ss)^{\frac{m}{2}} = \om^{\frac{m}{2}}$.

On a donc $\Ass^{(m)}_{\la,\mu}(\k) = (\la,\mu) \bullet \Ass^{(m)}_1(\k)$
d\`es que $\la \neq 0$ et $\mu \neq 0$.
Les cas
particuliers se ram\`enent donc \`a $\Phi \in \Ass^{(m)}_{0,1}(\k)$
et $\Phi \in \Ass^{(m)}_{1,0}(\k)$. Dans le premier cas
$\ss = s$ donc $\ss^2 = 1$ et $\tt = \Phi \om s e^{t_1} \Phi^{-1}$,
dans le deuxième $\ss = s e^{t_0}$ et $\tt = \Phi \om s \Phi^{-1}$
donc $\tt^2 = 1$.

Soit $\varphi$ l'automorphisme involutif de $W$ d\'efini par
$\varphi(s) = \om s$,
$\varphi(\om) = \om^{-1}$ et $f$ l'automorphisme de $\mathfrak{g}[\k]$,
donc de $\mathcal{A}(\k)$, d\'efini sur les
g\'en\'erateurs par $f(\t_r) = \t_{-r+1}$. L'automorphisme $f$
est $\varphi$-\'equivariant, au sens o\`u $f(g \bullet x) = \varphi(g)
\bullet f(x)$ pour tous $g \in W$, $x \in \mathsf{U} \mathfrak{g}(\k)$.
On en d\'eduit un automorphisme involutif $\J$ de $\mathfrak{B}(\k)$. Si $m$
est impair, il s'agit simplement de la conjugaison par $s
\om^{\frac{m-1}{2}}$. Si $m$ est pair en revanche, $\varphi$ est
un automorphisme ext\'erieur de $W$, donc il s'agit d'un
automorphisme ext\'erieur de $\mathfrak{B}(\k)$ qui \'etend
$\varphi$, et \'echange $\t_0$ et $\t_1$.

\subsection{Morphismes de monodromie}
On rappelle que l'on plonge $\mathfrak{B}'(\k)$ dans $\mathfrak{B}(\k)$
par $t_r \mapsto \t_r - T/m$.

Fixons $\la \in \k$, $\Phi \in \Ass^{(m)}_{\la}(\k) \subset
\mathcal{A}'(\k)$ et notons $\ss = s e^{\la t_0}$,
$\tt = \Phi \om s e^{\la t_1} \Phi^{-1} \in \mathfrak{B}'(\k)$,
$\sb = s e^{\la \t_0} \in \mathfrak{B}(\k)$,
$\tb = \Phi \om s e^{\la \t_1} \Phi^{-1}
\in \mathfrak{B}(\k)$. On a $\ss = \sb e^{-\la T/m}$,
$\tt = \tb e^{-T \la/m}$.

Si $m$ est impair, on d\'eduit du lemme \ref{relationssttimp} que
$$
\sb (\tb \sb)^{\frac{m-1}{2}} = (\tb \sb)^{\frac{m-1}{2}} \tb
= s \om^{\frac{m-1}{2}} \Phi^{-1} e^{\la T}
$$ 
et $(\sb \tb)^m = (\tb \sb)^m = e^{2 \la T}$. Si $m$
est pair, on d\'eduit du lemme \ref{relationssttpair} que
$(\sb \tb)^{\frac{m}{2}} = (\tb \sb)^{\frac{m}{2}}
= \om^{\frac{m}{2}} e^{ \la T}$ et $(\sb \tb)^m = e^{2 \la T}$.
Dans les deux cas, on a associ\'e \`a $\Phi \in \Ass^{(m)}_{\la}(\k)$
un morphisme $\widetilde{\Phi} : B \to \mathfrak{B}(\k)^{\times}$,
groupe form\'e des \'el\'ements grouplike de $\mathfrak{B}(\k)$,
d\'efini par $\sigma \mapsto \sb$, $\tau \mapsto \tb$. Il est clair
qu'il vérifie la condition fondamentale de l'introduction.

Remarquons que de la donn\'ee de $\PPhi$ on d\'eduit $\la$ par
la relation $\PPhi( (\sb \tb)^m) = e^{2 \la T}$. Si $m$
est impair on peut \'egalement en d\'eduire $\Phi$, par
$\PPhi( \sb (\tb \sb)^{\frac{m-1}{2}} ) = s \om^{\frac{m-1}{2}}
\Phi^{-1} e^{\la T}$. En revanche si $m$ est pair
cela n'est pas possible puisque $\Phi$ et $\Phi e^{t_1}$ donnent
lieu au m\^eme morphisme $\PPhi$.

On note $p : \mathfrak{B}(\k) \to \k W$ la projection
sur sa composante homogèene de degré 0. Comme $p \circ \PPhi(\sigma) = s$ et
$p \circ \PPhi(\tau) = \om s$ on en d\'eduit que $p \circ \PPhi = \pi$
et que
$\PPhi$
se restreint en un morphisme $\PPhi : P \to \exp \mathfrak{g}(\k)
= \mathcal{A}(\k)^{\times}$
qui s'étend en $\PPhi : P(\k) \to \exp \mathfrak{g}(\k)$
et finalement $\PPhi : B(\k) \to \mathfrak{B}(\k)$.
Il est clair que, si $f \in P'(\k)$, alors $\PPhi(f) \in \exp 
\mathfrak{g}'(\k)$. 

Nous établissons ici les propriétés élémentaires de ces morphismes.

\begin{prop} Soient $\la \neq 0$ et $\Phi \in \Ass^{(m)}_{\la}(\k)$. Le
morphisme $\PPhi : B(\k) \to \mathfrak{B}(\k)^{\times}$ est un isomorphisme
qui envoie $P(\k)$ (resp. $P'(\k)$) sur $\mathcal{A}(\k)^{\times}$
(resp. $\mathcal{A}'(\k)^{\times}$).
\end{prop}
\begin{proof}
On montre d'abord que $\PPhi : P(\k)
\to \exp \mathfrak{g}(\k)$ est un
isomorphisme. Comme c'est un morphisme de $\k$-sch\'emas en groupes
affines connexes (car pro-unipotents) il suffit de montrer que
l'application tangente $\dd \PPhi : L(\k) \to \mathfrak{g}(\k)$
est un isomorphisme de $\k$-alg\`ebres de Lie filtr\'ees compl\`etes.
Notons $u_r = \exp v_r$, $v_r \in L(\k)$. On a $\dd \PPhi (v_r) = 2 \la
t_r + w_r$, avec $w_r$ de valuation au moins 2. De plus,
$\PPhi(Z) = \exp (2 \la T)$, donc $\dd \PPhi (z) = 2 \la T$. D'autre
part, $u_{m-1} \dots u_1 u_0 = Z$ donc $v_0 + v_1 + \dots + v_{m-1} = z$.
On en d\'eduit un morphisme $\varphi : \mathfrak{g}(\k)
\to L(\k)$ tel que $\varphi(t_r) = v_r /2\la$. Alors $c = \dd \PPhi
\circ \varphi$ est un endomorphisme de $\mathfrak{g}(\k)$
tel que $c(t_r) = t_r + w_r$. On en d\'eduit que,
pour tout $x \in \mathfrak{g}(\k)$ la valuation de $c(x) -x $
est strictement sup\'erieure \`a la valuation de $x$ donc $c$ est
inversible et $\dd \PPhi : L(\k) \to \mathfrak{g}(\k)$
est un isomorphisme.

On a $p(\mathfrak{B}(\k)^{\times}) = W$ et $p \circ \PPhi = \pi$, donc
d\`es que $\la \neq 0$ on en
d\'eduit que $\PPhi : B(\k) \to  \mathfrak{B}(\k))^{\times}
$ est injectif. De plus, si $x \in \mathfrak{B}(\k)^{\times}$ et
$x_0 = p(x)$ il existe $\sigma \in B$ tel que $\pi(\sigma) = x_0$. On a
alors $p(\PPhi(\sigma)) = x_0$ d'où $\PPhi(\sigma^{-1}) x \in
\mathcal{A}(\k)^{\times}$ vaut $\PPhi(g)$ pour un certain $g \in P(\k)$
et $x = \PPhi(\sigma g)$, donc $\PPhi$ est surjectif et donc un
isomorphisme.

\end{proof}

\begin{lemme} \label{caractassimp} Soit $m \geq 3$ impair, $\la \in \k$,
$\Phi \in \exp \mathfrak{g}'(\k)$. On a $\Phi \in \Ass^{(m)}_{\la}(\k)$
si et seulement si il existe un morphisme $\PPhi : B(\k) \to
(\mathfrak{B}(\k))^{\times}$ (n\'ecessairement unique) tel que
$\PPhi(\sigma) = s e^{\la t_0}$ et
$$
\PPhi(\tau) = \Phi
\om s
e^{\la \t_1} \Phi^{-1}, \ \PPhi(O) = s \om^{\frac{m-1}{2}}
\Phi^{-1} e^{\la T},\ \Phi(O^2) = e^{2\la T}
$$
\end{lemme}
\begin{proof}
Si $\Phi \in \Ass^{(m)}_{\la}(\k)$, alors $\PPhi$ est bien d\'efini et
satisfait ces propri\'et\'es comme on vient de le montrer. Inversement,
posons $\PP = \Phi e^{\la t_1} \in \mathcal{A}'(\k)$. De
$\PPhi(O^2) = \PPhi(O)^2$ on d\'eduit (\ref{impeq2}). On montre alors
par calcul direct (cf. lemme \ref{relationssttimp}) que
$$
\PPhi\left(  (\tau \sigma)^{\frac{m-1}{2}} \right)
= \om^{\frac{m-1}{2}} \left( \om^{\frac{m+1}{2}} \bullet
\left( \prod_{q=0}^{\frac{m-3}{2}} \om^q \bullet \xi \right)
\right) e^{\la \frac{T}{m} (m-1)}
$$
avec $\xi = \PP (\om^{\frac{m+1}{2}} \bullet \PP)$ soit,
comme $O = \sigma (\tau \sigma)^{\frac{m-1}{2}}$,
$$
\PPhi(O) = s e^{\la t_0} \PPhi\left( (\tau \sigma)^{\frac{m-1}{2}}
\right) = s \om^{\frac{m-1}{2}} \Phi^{-1} e^{\la T}.
$$
On en d\'eduit  
$$
\om^{\frac{m+1}{2}} \bullet
\left( \prod_{q=0}^{\frac{m-3}{2}} \om^q \bullet \xi \right)
= e^{-\la \t_1} \Phi^{-1} = \PP^{-1}
$$
d'o\`u (\ref{impeq3}).
\end{proof}
\begin{lemme} \label{caractasspair}
Soit $m \geq 4$ pair, $\la \in \k$,
$\Phi \in \exp \mathfrak{g}'(\k)$. On a $\Phi \in \Ass^{(m)}_{\la}(\k)$
si et seulement si il existe un morphisme $\PPhi : B(\k) \to
\mathfrak{B}(\k)^{\times}$ (n\'ecessairement unique) tel que
$\PPhi(\sigma) = s e^{\la \t_0}$ et
$$
\PPhi(\tau) = \Phi
\om s
e^{\la \t_1} \Phi^{-1}, \ \PPhi(O) = \om^{\frac{m}{2}} e^{\la T},
\ \Phi(O^2) = e^{2\la T}
$$
\end{lemme}
\begin{proof} Si $\Phi \in \Ass^{(m)}_{\la}(\k)$, $\PPhi$ v\'erifie bien
ces propri\'et\'es. Inversement, posons $\xi = \Phi e^{\la t_1} 
(\om s \bullet \Phi^{-1}) e^{\la t_2}$. On montre par calcul direct (cf.
lemme \ref{relationssttpair}) que
$$
\PPhi \left( ( \tau \sigma)^{\frac{m}{2}} \right)
= \left( \prod_{q=0}^{\frac{m}{2} -1} \om^q \bullet  \xi \right)
\om^{\frac{m}{2}} e^{\la T}
$$
d'o\`u (\ref{paireq3})
puisque $\PPhi( (\tau \sigma)^{\frac{m}{2}}) = \PPhi(O) = \om^{\frac{m}{2}}
e^{\la T}$.
\end{proof}

Ces deux lemmes nous donnent une d\'emonstration de la
proposition \ref{assplusmoins} alternative au calcul direct. En effet,
soit $\Phi \in \Ass^{(m)}_{\la}(\k)$
et $\Gamma$ l'automorphisme de $B$ d\'efini par $\sigma \mapsto \sigma^{-1}$,
$\tau \mapsto \tau^{-1}$ (\og image miroir de la tresse \fg). Le
morphisme $\Gamma$ s'\'etend \`a $B(\k)$ et permet de d\'efinir
$\widetilde{\Psi} = \PPhi \circ \Gamma$. On a $\widetilde{\Psi}(\sigma)
= s e^{-\la t_0}$, $\widetilde{\Psi}(\tau) = \Phi \om s
e^{-\la \t_1} \Phi^{-1}$,
$\widetilde{\Psi}(O^2) = \PPhi(O^{-2}) = e^{-2\la T}$,
$\Gamma(O) = O^{-1}$ donc $\widetilde{\Psi}(O) = \om^{\frac{m}{2}}
e^{-\la T}$ si $m$ est pair,
$$
\widetilde{\Psi}(O) = e^{-\la T} \Phi^{-1} s \om^{\frac{m-1}{2}}
= s \om^{\frac{m-1}{2}}  (s \om^{\frac{m-1}{2}} \bullet \Phi)
e^{-\la T} = s \om^{\frac{m-1}{2}} \iota(\Phi)^{-1} e^{-\la T}
$$
si $m$ est impair. On en d\'eduit que $\widetilde{\Psi}$
est le morphisme associ\'e \`a $\Phi$ et au scalaire $-\la$, et
que $\Phi \in \Ass^{(m)}_{-\la}(\k)$.

\section{Groupe d'automorphismes de $B(\k)$}
\setcounter{equation}{0}

\subsection{Groupe sp\'ecial d'automorphismes}

On suppose $m \geq 3$. On note
$$O = \underbrace{\sigma \tau \sigma \dots}_{\mbox{$m$ facteurs}}.
$$
Si $m$ est impair, $\pi(O)$ est une r\'eflexion. Que $m$ soit pair
ou impair, $\pi(O)$ est d'ordre 2, donc $O$ est justifiable
de la notation $O^{<\la>}$. 

\begin{definition} On note $G^{(m)}(\k)$ l'ensemble des couples
$(\la, f) \in \kt \times P'(\k)$ tels qu'il existe un morphisme
$\tilde{f} : B \to B(\k)$ satisfaisant \`a
\begin{eqnarray} 
\sigma \mapsto \sigma^{<\la>} \label{gma} \\
\tau \mapsto f^{-1} \tau^{<\la>} f \label{gmb} \\
O \mapsto O^{<\la>} f \mbox{\ \ si $m$ est impair,\ \ } O \mapsto O^{<\la>}
 \mbox{\ \ si $m$ est pair.} \label{gmc} \\
O^2 \mapsto (O^2)^{\la} \label{gmd} 
\end{eqnarray}
Si $m$ est pair, (\ref{gmc}) implique naturellement (\ref{gmd}). On \'etend naturellement $\tilde{f}$ en un automorphisme continu pur
de $B(\k)$.
\end{definition}
L'automorphisme $\tilde{f}$ est bien d\'etermin\'e par $f$. Inversement,
si $m$ est impair, la condition (\ref{gmc}) montre que $f$ est d\'etermin\'e
par $\tilde{f}$. On en d\'eduit dans ce cas que l'application $f \mapsto
\tilde{f}$ de $G^{(m)}(\k)$ vers $Aut(B(\k))$ est injective.
Il n'en est pas de m\^eme si $m$ est pair, puisqu'alors $(\la, \tau^2 f)$
correspond au m\^eme morphisme $\tilde{f}$.

Soient maintenant $\tilde{f}_1$, $\tilde{f}_2$ les automorphismes
correspondant \`a deux \'el\'ements $(\la_1,f_1)$ et $(\la_2,f_2)$
de $G^{(m)}(\k)$ et $\tilde{f} = \tilde{f}_1 \circ \tilde{f}_2$.
On a
$
\tilde{f}(\sigma) =  \sigma^{<\la>}
$
avec $\la = \la_1 \la_2$, $\tilde{f}(\tau) = f^{-1} \tau^{<\la>} f$
avec $f = f_1 \tilde{f}_1(f_2)$ et $\tilde{f}(O)$ vaut
$O^{<\la>} f$ si $m$ est impair,  $O^{<\la>}$ si $m$ est pair.
On en d\'eduit une loi de composition interne
associative sur $G^{(m)}(\k)$ donn\'ee par la formule
$$
(\la_1, f_1) * (\la_2 ,f_2) = (\la_1 \la_2, f_1 \tilde{f}_1(f_2)).
$$
Comme de plus $(1,1)$ est un \'el\'ement neutre \'evident et que 
$(\la, f) * (\la^{-1}, (\tilde{f}^{-1})(f^{-1})) = (1,1)$
on vient de munir $G^{(m)}$ d'une loi de groupe, l'identifiant
\`a un sous-groupe des automorphismes de $B(\k)$.

Pour simplifier les notations, on notera \`a $\la$ fix\'e
$\sss = (\sigma^2)^{\frac{\la-1}{2}} = \sigma^{-1} \sigma^{<\la>} \in P(\k)$
et $\ttt = \tau^{-1} \tau^{<\la>}$, $\ooo = O^{-1}O^{<\la>}$.

On peut traduire la définition de $G^{(m)}(\k)$ par des conditions
al\-gé\-bri\-ques sur les couples $(\la,f)$. Si $m$ est impair les conditions
que l'on obtient ici peuvent se simplifier (voir \ref{ssalglie}).
\begin{lemme}
$G^{(m)}(\k)$ est l'ensemble des couples $(\la,f) \in \kt \times P'(\k)$
tels que
\begin{eqnarray}
O \bullet f = f^{-1} \mbox{\ \   si $m$ est impair.} \label{gexp1} \\
\left[ \prod_{k=\frac{m-1}{2}}^{1} (\tau \sigma)^{-k} \bullet L
\right] \sss = \ooo f \mbox{\ \  si $m$ est impair,} \label{gexp2} \\
\nonumber \prod_{k=\frac{m}{2}-1}^0 ( (\sigma \tau)^{-k} \bullet M) = \ooo
\mbox{\ \   si $m$ est pair.}
\end{eqnarray}
avec $L = \sss f^{-1} \ttt (\tau \bullet f)$, $M = (\tau^{-1} \bullet
(\sss f^{-1})) \ttt f$.
\end{lemme}
\begin{proof}
Soit $(\la,f) \in G^{(m)}(\k)$. Si $m$ est impair, les conditions (\ref{gmc})
et (\ref{gmd}) impliquent que $
(O^2)^{\la} = \tilde{f}(O)^2$ vaut
$$
(O^{<\la>} f)^2 = (O^2)^{\frac{\la-1}{2}} O
f O^{-1} O^2 (O^2)^{\frac{\la-1}{2}}= (O^2)^{\la} (O \bullet f) f
$$
car $O^2$ est central, d'o\`u (\ref{gexp1}). Soit $V_r = \sigma \tau \dots $
($r$ facteurs) et $E_r \in P(\k)$ tel que $\tilde{f}(V_r) = V_r E_r$.
On a $V_1 = \sigma$, $\tilde{f}(V_1) = \sigma^{<\la>} = \sigma \sss
= V_1 \sss$ d'o\`u $E_1 = \sss$. On a 
$E_{2r+1} = (\sigma^{-1} \bullet E_{2r}) \sss$ et 
$E_{2r+2} = (\tau^{-1} \bullet E_{2r+1}f^{-1}) 
\ttt f$. On en déduit
$E_{2r+2} = ((\tau^{-1} \sigma^{-1}) \bullet E_{2r}) M$ et
$F_{2r+3} = (\tau \sigma)^{-1} \bullet F_{2r+1} L$ où l'on a posé
$E_{2r+1} = F_{2r+1} \sss$,
et enfin que
$$
E_m = \prod_{k = \frac{m}{2} -1}^0 \left( (\sigma \tau)^{-k} \bullet M
\right) \mbox{ ou } E_m = \left[ \prod_{k = \frac{m-1}{2}}^1  (\tau \sigma)^{-k} \bullet L)  
\right] \sss
$$
suivant que $m$ est pair ou impair.
Comme $V_m = O$,
on d\'eduit alors (\ref{gexp2}) de (\ref{gmc}). Inversement, on montre
par la m\^eme m\'ethode que, si $(\la,f)$ v\'erifie (\ref{gexp1})
et (\ref{gexp2}), alors $(\la,f) \in G^{(m)}(\k)$.
\end{proof}

\subsection{Action de $G^{(m)}(\k)$ sur $\Ass^{(m)}(\k)$}

Soient $\Phi \in \Ass^{(m)}_{\mu}(\k)$ et $(\la,f) \in G^{(m)}(\k)$. Leur sont
associ\'es des morphismes $\PPhi$ et $\tilde{f}$ :
$$
B(\k) \stackrel{\tilde{f}}{\to} B(\k) \stackrel{\PPhi}{\to} 
\mathfrak{B}(\k) ^{\times}.
$$
On a
$$
(\PPhi \circ \tilde{f})(\sigma) = \PPhi(\sigma^{<\la>})
= \PPhi(\sigma) \PPhi( (\sigma^2)^{\frac{\la-1}{2}} )
= s e^{\mu t_0} \left( e^{2\mu t_0} \right)^{\frac{\la -1}{2}} = s e^{\la \mu
t_0}
$$
et
$$
(\PPhi \circ \tilde{f} )(\tau) = \PPhi( f^{-1} \tau^{<\la>}
f)
= \PPhi(f)^{-1} \Phi \om s e^{\la \mu t_1}
\Phi^{-1} \PPhi(f).
$$
On pose alors $\Phi .f =  \PPhi(f)^{-1} \Phi$.
On a donc
$\Phi . f
\in \exp  \mathfrak{g}'(\k)$.

\begin{lemme} \label{actiongmass} Si $\Phi \in \Ass^{(m)}_{ \mu}(\k)$, $(\la, f) \in G^{(m)}(\k)$
alors $\Phi .f \in \Ass^{(m)}_{\la \mu}(\k)$.
\end{lemme}
\begin{proof}
D'apr\`es ce qui pr\'ec\`ede et les lemmes \ref{caractassimp} et
\ref{caractasspair} 
il suffit de calculer les images par $\PPhi
\circ
\tilde{f}$ de $O$ et $O^2$.
On a $\PPhi \circ \tilde{f} (O^2) = \PPhi((O^2)^{\la})=\PPhi(O^2)^{\la} = e^{2 \la \mu T}$.
Si $m$ est impair,
$$
\begin{array}{lcl}
\PPhi \circ \tilde{f} (O) & = & \PPhi( O^{<\la>} f) = \PPhi(O)
\PPhi(O^2)^{\frac{\la-1}{2}} \PPhi(f) \\
& = & s \om^{\frac{m-1}{2}} \iota(\Phi)^{-1} e^{\mu T} e^{\mu (\la -1) T}
\PPhi(f) \\
& = & s \om^{\frac{m-1}{2}} \iota( \Phi .f)^{-1}  e^{\la \mu T}.
\end{array}
$$
Si $m$ est pair, $\PPhi \circ \tilde{f}(O) = \PPhi (O^{<\la>}) = 
\PPhi (O) e^{\mu T(\la-1)} = \om^{\frac{m}{2}} e^{\la \mu T}$.
\end{proof}

Si $\Phi \in \Ass^{(m)}_{\mu}(\k)$ et $f \in P'(\k)$, on peut
toujours
d\'efinir $\Phi . f =  \PPhi(f)^{-1} \Phi \in \exp \mathfrak{g}'(\k)$.

\begin{lemme} \label{lemmetrans}
Soit $\Phi \in \Ass^{(m)}_{\mu}(\k)$  et $f \in P'(\k)$. Si
$\Phi .f
\in \Ass^{(m)}_{\la \mu}(\k)$ avec $\la \mu \neq 0$ alors $(\la,f) \in G^{(m)}(\k)$.
\end{lemme}
\begin{proof}
D\'efinissons $\tilde{f}(\sigma) = \sigma^{<\la>}$,
$\tilde{f}(\tau) = f^{-1} \tau^{<\la>} f$. Cette application
s'\'etend en un homomorphisme de groupes du groupe libre sur $\sigma,
\tau$ vers $B(\k)$, telle que $\PPhi \circ \tilde{f} = \widetilde{\Phi .f}$,
donc $\tilde{f}$ se factorise par $B$ parce que $\PPhi$ est injectif
sur $B(\k)$. Par d\'efinition, $\tilde{f}$ satisfait (\ref{gma}),
(\ref{gmb}), et $\PPhi \circ \tilde{f} = \widetilde{\Phi.f}$. On a
$$
\PPhi \circ \tilde{f} (O^2) = e^{2 \la \mu T} = \PPhi(O^2)^{\la}
= \PPhi( (O^2)^{\la})
$$
d'o\`u $\tilde{f}(O^2) = (O^2)^{\la}$, soit (\ref{gmd}), par
injectivit\'e de $\PPhi$. De m\^eme, si $m$ est impair,
$$
\PPhi \circ \tilde{f}(O) = s \om^{\frac{m-1}{2}} \Phi^{-1}
\PPhi(f) e^{\la \mu T} = s \om^{\frac{m-1}{2}} \Phi^{-1}
e^{\la \mu T} \PPhi(f) = \PPhi(O^{<\la>} f)$$
d'o\`u $\tilde{f}(O) = O^{<\la>} f$ par injectivit\'e de $\PPhi$,
soit (\ref{gmc}). Si $m$ est pair, on conclut de m\^eme
\`a partir de $\PPhi \circ \tilde{f}(O) = \PPhi(O^{<\la>}) = O^{<\la \mu>}$.
\end{proof}

On note $\Ass^{(m)}_*(\k)$ la r\'eunion des couples $(\la , \Phi)$
pour $\la \in \kt$, $\Phi \in \Ass^{(m)}_{\la}(\k)$. D'apr\`es le lemme
\ref{actiongmass} on a une action \`a droite de $G^{(m)}(\k)$
sur $\Ass^{(m)}_*(\k)$ par $(\mu,\Phi).(\la,f) = (\mu \la, \Phi.f)$.

\begin{prop} \label{freetransitive}
L'action de $G^{(m)}(\k)$ sur $\Ass^{(m)}_*(\k)$ est libre et transitive.
\end{prop}
\begin{proof}
On d\'emontre d'abord la transitivit\'e. Soient $\Phi \in \Ass^{(m)}_{\la}(\k)$, $\Phi' \in \Ass^{(m)}_{\mu}(\k)$ avec
$\la \mu \neq 0$. On a $\Phi' \Phi^{-1} \in \exp
\mathfrak{g}'(\k)$. Posons $f = \PPhi^{-1}(\Phi' \Phi^{-1}) \in P'(\k)$.
On a alors $\Phi' = \PPhi(f)^{-1} \Phi$
donc $\Phi' = \Phi.f$ et $(\mu/\la, f) \in G^{(m)}(\k)$ d'apr\`es le lemme
\ref{lemmetrans}, donc l'action est transitive. Si $\Phi = \Phi.f$
avec $(\la, f) \in G^{(m)}(\k)$, alors $\PPhi = \PPhi \circ 
\tilde{f}$ donc $\tilde{f}$ est l'identité par injectivit\'e de $\PPhi$. Cela
implique $\la = 1$ et, au moins si $m$ est impair, $f = 1$. En g\'en\'eral,
comme $f \in P'(\k)$ on a $\Phi .f = \PPhi(f^{-1}) \Phi$
donc $\Phi .f = \Phi \Rightarrow \PPhi(f) = 1 \Rightarrow f=1$
puisque $\PPhi$ est injectif ($\la = 1 \neq 0$).
\end{proof}

\section{Existence d'un associateur rationnel} \label{section6}
\setcounter{equation}{0}

\subsection{G\'en\'eralit\'es}

On fixe $m \geq 3$, ce qui nous permet d'\^oter
l'exposant $(m)$ des notations.
L'action de $\kt$ sur $\mathfrak{B}'(\k)$
qui \`a $\alpha \in \kt$ et $t_i \in \mathfrak{B}'(\k)$ associe
$\alpha t_i$ induit une action de $\kt$ sur $\Ass_*(\k)$ : \`a
$\alpha \in \kt$ et $(\la, \Phi) \in \Ass_*(\k)$ on associe $(\la \alpha,
\Phi_{\alpha})$ o\`u $\Phi_{\alpha}$ est l'image de $\Phi \in
\mathfrak{B}'(\k)$ par l'action de $\alpha$. Le quotient $\Ass_*(\k)/
\kt$ s'identifie \`a $\Ass_1(\k)$. On note
$j : \Ass_1(\k) \hookrightarrow \Ass_*(\k)$ et $q : \Ass_*(\k)
\onto \Ass_1(\k)$ les injections et surjections canoniques.
Si $\Phi \in \Ass_1(\k)$ et $f \in G(\k)$ on pose
$\Phi \star f = q(\Phi.f)$. C'est une action de $G(\k)$ sur
$\Ass_1(\k)$.

Supposons $\Ass_1(\k) \neq \emptyset$ et choisissons
$\Phi \in \Ass_1(\k)$. Consid\'erons la suite
\begin{eqnarray} \label{G1GK}
1 \to G_1(\k) \to G(\k) \stackrel{\nu}{\to} \kt \to 1
\end{eqnarray}
o\`u $\nu(\la,f) = \la$. Elle est exacte si et seulement si $\nu$
est surjective. Pour tout $\la \in \kt$ on a $\Phi_{\la}
\in \Ass_{\la}(\k) \subset \Ass_*(\k)$ donc il existe un unique
$g = (\la, \check{g}) \in G(\k)$ tel que $\Phi_{\la} = 
\Phi. \check{g}$. Posant $\Theta_{\Phi}(\la) = g$ on a donc $\Theta_{\Phi}
(1) = (1,1)$ et $\nu \circ \Theta_{\Phi}(\la) = \la$, donc
(\ref{G1GK}) est exacte et \`a chaque $\Phi \in \Ass_1(\k)$ correspond
une section de $\nu$. On a
$$
\begin{array}{lcl}
\Theta_{\Phi}(\kt) & = & \{ g \in GT(\k) \ \mid \ \exists \la \in \kt \ \ \Phi_{\la} = \Phi.g \} \\
& = & \{ g \in GT(\k) \ \mid \ q(\Phi.g) = \Phi \} \\
& = & \{ g \in GT(\k) \ \mid \ \Phi = \Phi \star g \}
\end{array}
$$
ainsi $\Theta_{\Phi}(\kt)$ est le stabilisateur de $\Phi$ pour
l'action de $GT(\k)$ sur $\Ass_1(\k)$. D'autre part, notant
$\mathcal{G}(\k)$ et $\mathcal{G}_0(\k)$ les alg\`ebres de Lie
respectives de $G(\k)$ et $G_1(\k)$, on en d\'eduit une suite
exacte
\begin{eqnarray} \label{g1gk}
0 \to \mathcal{G}_0(\k) \to \mathcal{G}(\k) \stackrel{\dd \nu}{\to} \k
\to 0
\end{eqnarray}
et une section $\dd \Theta_{\Phi}$ de $\dd \nu$. Notant $(\la, g_{\la})
= \Theta_{\Phi}(\la)$, on a
\begin{eqnarray} \label{eqstabl}
\Phi_{\la} = \Phi.g_{\la} = \PPhi(g_{\la}^{-1}) \Phi = (\PPhi \circ
\Theta_{\Phi}(\la^{-1})) \Phi
\end{eqnarray}

D'autre part $G(\k) \subset \kt \times P'(\k)$, d'alg\`ebre de Lie
$\k \times L'(\k)$. Posons $\dd \Theta_{\Phi} (1) = (1,\psi) \in
\k \times L'(\k)$, d'o\`u $\dd \Theta_{\Phi}(l) = (l,l\psi)$. Pour
tout $\Phi \in \mathcal{A}'(\k)$, consid\'erons le 
morphisme de sch\'emas $\mathbbm{G}_m \to \mathcal{A}'(.)$ qui \`a
$\alpha \in \kt $ associe $\Phi_{\alpha}\in \mathcal{A}'(\k)$. Sa diff\'erentielle
permet de construire une d\'erivation $\partial$ de $\mathcal{A}'(\k)$
d\'efinie par $\Phi_{1 + \eps l} = \Phi + \eps l \partial \Phi$
pour $\eps^2 = 0$ et $l \in \k$. Diff\'erenciant (\ref{eqstabl})
en $\la = 1$ on obtient alors $\partial \Phi = - \left( \dd \PPhi
\circ \dd \Theta_{\Phi} (1) \right) \Phi$ soit
\begin{eqnarray} \label{maineq}
 \left( \partial \Phi \right) \Phi^{-1} = - \dd \PPhi(\psi)
\end{eqnarray}
L'existence d'un associateur transcendant implique $\Ass_1(\C) \neq 
\emptyset$, donc les suites (\ref{G1GK}) et (\ref{g1gk}) sont
exactes, c'est-\`a-dire que $\dd \nu$ est non nulle
sur les points rationnels complexes. Comme $\dd \nu$ est d\'efinie
sur $\Q$, il s'ensuit que $\dd \nu$ est non nulle sur $\Q$, donc
que (\ref{g1gk}) est exacte pour $\k = \Q$. Il existe donc
$(1,\psi) \in \mathcal{G}(\Q)$ tel que $\psi \neq 0$.

Dans le cas pair, nous aurons besoin d'introduire un
sous-groupe propre de $G(\k)$.
On note $J$ l'automorphisme de $B$ d\'efini par $J(\sigma) = \tau$,
$J(\tau) = \sigma$. Quand $m$ est impair, il s'agit de la
conjugaison par $O$. Lorsque $m$ est pair, il s'agit d'un
automorphisme \emph{ext\'erieur} de $B$. Dans les deux cas il
se restreint en un automorphisme ext\'erieur de $P$, et il
se prolonge en un automorphisme de $B(\k)$ qui laisse stable
$P(\k)$ et $P'(\k)$.

On note $\J$ l'automorphisme de $\mathfrak{B}(\k)$ d\'efini par
$\J(s) = \om s$, $\J(\om) = \om^{-1}$, $\J(t_{r}) = r_{-r+1}$. Quand
$m$ est impair il s'agit de la conjugaison par $s \om^{\frac{m-1}{2}}$,
quand $m$ est pair d'un automorphisme involutif ext\'erieur de
$\mathfrak{B}(\k)$.

On d\'efinit
$$
\begin{array}{l}
\Ass'_{\la}(\k) = \{ \Phi \in P'(\k) \ \mid \ (\la, \Phi) \in \Ass_{\la}(\k),
\ \ \ \J(\Phi) = \Phi^{-1} \} \\
\Ass'(\k) = \{ (\la,\Phi) \in \Ass(\k) \ \mid \ \J(\Phi) = \Phi^{-1} \} \\
G'(\k) = \{ (\la,f) \in G(\k) \ \mid \ J(f) = f^{-1} \}
\end{array}
$$
et $\Ass'_*(\k) = \Ass'(\k) \cap \Ass_*(\k)$. Soit $\Phi \in \Ass'_{\la}
(\k)$ et $\PPhi : B(\k) \to \mathfrak{B}(\k)^{\times}$
le morphisme associ\'e. On a $\J(\PPhi(\sigma)) = \om s e^{\la t_1}
= \Phi^{-1} \bullet \PPhi(\tau)$ et $\J(\PPhi(\tau)) = \Phi^{-1} \bullet
\PPhi(\sigma)$. On en d\'eduit
\begin{eqnarray}
\label{commJphi} \J \circ \PPhi = \Ad(\Phi^{-1}) \circ \PPhi \circ J
\end{eqnarray}
Si $(\la,f) \in G'(\k)$ et $\tilde{f} : B \to B(\k)$ d\'esigne
le morphisme
associ\'e, on a 
\begin{eqnarray}
\label{commJf} \tilde{f} \circ J = \Ad(f^{-1}) \circ J \circ \tilde{f}
\end{eqnarray}
On en d\'eduit que $G'(\k)$ est un sous-groupe de $G(\k)$. De plus,
on d\'eduit de (\ref{commJphi}) que $\Ass'(\k)$ est stable sous
l'action de $G'(\k)$. Enfin, si $\Phi \in \Ass'_{\mu}(\k)$ avec
$\mu \neq 0$ et $f \in P'(\k)$ tel que $\Phi .f \in \Ass'_{\la \mu}(\k)$,
alors on sait d'apr\`es le lemme \ref{lemmetrans} que $(\la,f)
\in G(\k)$. On a donc $\J(\Phi) = \Phi^{-1}$ et $\J \circ \PPhi(f^{-1})
\J(\Phi) = \Phi^{-1} \PPhi(f)$. On en d\'eduit
$\J \circ \PPhi(f^{-1}) = \Ad(\Phi^{-1}) \circ \PPhi \circ J(f^{-1})$,
d'o\`u $f = J(f)^{-1}$ d'apr\`es (\ref{commJphi}) et l'injectivit\'e de
$\PPhi$. Ainsi l'action de $G'(\k)$ sur $\Ass'_*(\k)$ est libre et
transitive. Si $m$ est impair on a évidemment $\Ass(\k) = \Ass'(\k)$,
$G(\k) = G'(\k)$, etc.

\begin{lemme} \label{lemmtranspair} L'associateur transcendant $\Phi_0$ appartient \`a $\Ass'_{\ii
\pi}(\C)$.
\end{lemme}
\begin{proof} On note $\j$ le g\'en\'erateur de $\Z/2\Z$. On a des
actions de $\Z/2\Z$ sur $\mathcal{A}'(\k)$ et $\R^2$, en faisant
agir $\j$ sur $\mathcal{A}'(\k)$ par $\J$, et sur $\R^2$
comme r\'eflexion par rapport \`a la droite $\Delta$ d'\'equation
$y =  -x \, \cotg \frac{\theta}{2}$
(cf. figure \ref{exterieur}).
\begin{figure}
\begin{center}
\includegraphics{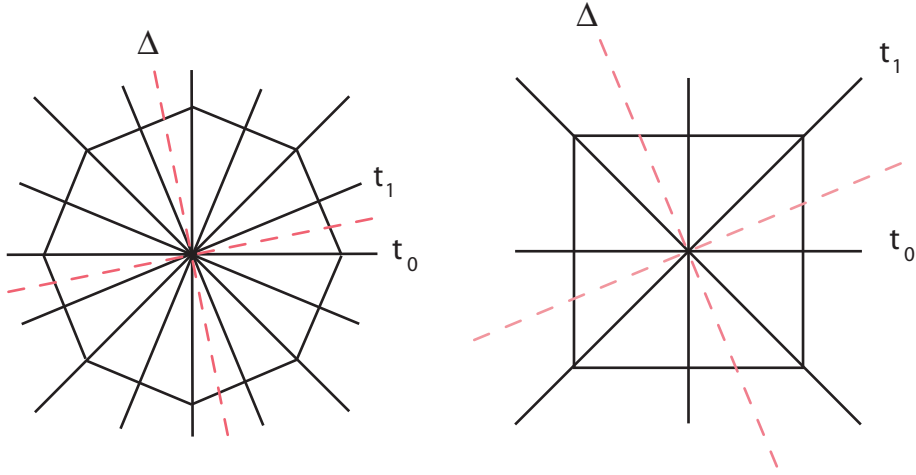}
\end{center}
\caption{L'automorphisme ext\'erieur $\j$ pour $m=8$ et $m=4$} \label{exterieur}
\end{figure}

La forme diff\'erentielle
$\Omega'$ est alors invariante par $\j$, et si $F_{0,+}$, $F_{0,-}$
sont d\'efinis
comme en section 3 on a $\j.F_{0,+} = F_{1,-}$ par comparaison
des comportements asymptotiques. On d\'eduit alors de $F_{1,-} = F_{0,+}
\Phi_0$ que $\J(\Phi_0) = \Phi_0^{-1}$.
\end{proof}

Le raisonnement précédent s'applique alors, en rempla\c cant
$G(\k)$, $G_1(\k)$ par $G'(\k)$, $G'_1(\k) = G'(\k) \cap G_1(\k)$
et $\mathcal{G}(\k)$, $\mathcal{G}_0(\k)$ par les alg\`ebres de Lie
$\mathcal{G}'(\k)$, $\mathcal{G}'_0(\k)$ de $G'(\k)$, $G'_1(\k)$.

On cherche donc \`a r\'esoudre (\ref{maineq}) avec $(1,\psi) \in
\mathcal{G}'(\k)$. Pour ce faire, il nous faut d'abord
définir $\PPhi$ pour tout $\Phi \in \mathcal{A}'(\k)$ de terme
constant égal à 1.

\subsection{Le morphisme $\PPhi$}

On rappelle que, si $A$ est une alg\`ebre de Hopf (compl\'et\'ee),
$A^{\times}$ d\'esigne le groupe des \'el\'ements grouplike
de $A$, et que si $A$ est une alg\`ebre $A^*$ d\'esigne le
groupe des \'el\'ements inversibles de $A$. 

On suppose $m$ impair. On a d\'efini des \'el\'ements $u_{2r} = 
(\tau \sigma)^r \bullet \sigma^2$, $u_{2r+1} = (\tau \sigma)^r \bullet
\tau^2$ et not\'e que $P'(\k)$ est engendr\'e par $\tilde{u}_r
= u_r Z^{\frac{-1}{m}}$ pour $0 \leq r \leq m-1$ soumis \`a l'unique
relation $\tilde{u}_{m-1} \dots \tilde{u}_0 = 1$. Il s'identifie
donc \`a la compl\'etion $\k$-prounipotente du groupe libre
sur toute partie \`a $m-1$ \'el\'ements de $\{ \tilde{u}_0,\dots,
\tilde{u}_{m-1} \}$. D'autre part, l'automorphisme $J$ laisse
stable $P'(\k)$. On
montre facilement les identit\'es suivantes
\begin{eqnarray}
\forall r \in [0,\frac{m}{2} ] \ \ \ u_{m-2r} = (\tau \sigma)^{-r}
 \bullet \sigma^2 \\
\forall r \in [0,\frac{m}{2} ] \ \ \ u_{m-2r+1} = (\tau \sigma)^{-r} \bullet \tau^2
\end{eqnarray}

On note $x_1 = u_1$, $x_0 = u_0$, $x_r = u_r \dots u_2 u_1$
pour $r \geq 1$, $x_r = u_0 u_{m-1} \dots u_{m+r}$ pour $r \leq 0$. On
montre par r\'ecurrence sur $r$ que
\begin{eqnarray}
\label{descrxd} \forall r \geq 1 \ \ x_{2r} = \underbrace{\tau \sigma \tau \dots \tau}_{2r-1}
\sigma^2 \underbrace{\tau \sigma \tau \dots \tau}_{2r-1} = (\tau \sigma)^r
(\sigma \tau)^r \\
\forall r \geq 0 \ \ x_{2r+1} = ( \underbrace{\tau \sigma \tau \dots
\tau}_{2r+1} )^2 = ((\tau \sigma)^r \tau)^2 \\
\forall r \geq 0 \ \ x_{-2r} = ( \underbrace{\sigma \tau \sigma \dots
\sigma}_{2r+1} )^2 = ((\sigma \tau)^r \sigma)^2 \\
\label{descrxf}\forall r \geq 1 \ \ x_{-2r+1} = \underbrace{\sigma \tau \sigma \dots
\sigma}_{2r-1}
\tau^2 \underbrace{\sigma \tau \sigma \dots \sigma}_{2r-1} =(\sigma \tau)^r (\tau \sigma)^r
\end{eqnarray}
D'autre part, on a $u_1 = x_1$, $u_0 = x_0$, $u_r = x_r x_{r-1}^{-1}$
pour $2 \leq r \leq \frac{m}{2}$, $u_{m-r} = x_{1-r}^{-1} x_{-r}$
pour $1 \leq r < \frac{m-1}{2}$. Posons $\tilde{x}_r = x_r
Z^{\frac{-r}{m}}$ pour $r \geq 1$, $\tilde{x}_{-r} = x_{-r} Z^{-\frac{r+1}{m}}$
pour $r \leq 0$ de telle sorte que $\tilde{x}_r \in P'(\k)$.
On a alors
\begin{eqnarray}
\forall r \geq 1, \ \ J(\tilde{x}_r) = \tilde{x}_{-r+1}. 
\end{eqnarray}

Si $m$ est impair, consid\'erons les deux ensembles
$$ \{ \tilde{u}_0, \tilde{u}_1,\dots, \tilde{u}_{\frac{m-1}{2}},
\tilde{u}_{\frac{m+3}{2}} , \dots, \tilde{u}_{m-1}  \}
= \{ \tilde{u}_0, \dots ,\widehat{\tilde{u}_{\frac{m+1}{2}}},\dots,
\tilde{u}_{m-1} \}
$$
et $\{ \tilde{x}_{-\frac{m-3}{2}},\dots,\tilde{x}_{-1}, \tilde{x}_0,
\tilde{x}_1,\dots, \tilde{x}_{\frac{m-1}{2}} \}$. Le groupe $P'(\k)$
est libre sur le premier ensemble, donc sur le deuxi\`eme d'apr\`es
les identit\'es pr\'ec\'edentes. Pour construire un morphisme
$\PPhi : P'(\k) \to \mathcal{A}'(\k)$ il suffit donc de déterminer
ses valeurs sur l'un de ces deux ensembles.

Si $m$ est pair, la relation $\tilde{u}_{m-1} \dots \tilde{u}_0 = 1$
est équivalente à $\tilde{x}_{\frac{m}{2}}
\tilde{x}_{1 - \frac{m}{2}} = 1$. Pour construire un morphisme
$\PPhi : P'(\k) \to \mathcal{A}'(\k)$ il suffit donc de déterminer
des valeurs $\PPhi(\tilde{x}_r)$ telles que
$\PPhi(\tilde{x}_{\frac{m}{2}})
\PPhi(\tilde{x}_{1 - \frac{m}{2}}) = 1$.

Pour tout $\Phi \in \mathcal{A}'(\k)^*$ tel que $\Phi \equiv
1 \mod \mathcal{A}'_1(\k)$, on note
$$
\sigma_{\Phi} = s e^{t_0}, \ \ \tau_{\Phi}   = \Phi \om s e^{t_1}
\Phi^{-1} \in \mathfrak{B}'(\k)
$$
et on d\'efinit un morphisme $\PPhi : P'(\k)\to \mathcal{A}'(\k)^*$
par les formules de la table \ref{tablephix}. On en déduit les
valeurs de $\PPhi$ sur les éléments $\tilde{u}_r$ (table \ref{tablephiu}).
\begin{table}
$$
\begin{array}{lclcl}
\PPhi(\tilde{x}_{2r}) & = & (\tau_{\Phi} \sigma_{\Phi})^r (\sigma_{\Phi}
\tau_{\Phi})^r & \mbox{ si } & 1 \leq 2r < \frac{m}{2} \\
\PPhi(\tilde{x}_{2r+1}) & = &( (\tau_{\Phi} \sigma_{\Phi})^r \tau_{\Phi})^2
 & \mbox{ si } & 1 \leq 2r +1 < \frac{m}{2} \\
\PPhi(\tilde{x}_{-2r}) & = &( (\sigma_{\Phi} \tau_{\Phi})^r \sigma_{\Phi})^2
 & \mbox{ si } & 1-\frac{m}{2}  < -2r  \leq 0 \\
\PPhi(\tilde{x}_{-2r+1}) & = & (\sigma_{\Phi} \tau_{\Phi})^r (\tau_{\Phi}
\sigma_{\Phi})^r & \mbox{ si } & 1-\frac{m}{2}  < -2r+1  \leq 0 \\
\PPhi(\tilde{x}_{\frac{m}{2}}) & = & \om^{\frac{m}{2}}
(\sp^{-1} \tp^{-1})^{\frac{m}{4}} (\sp \tp)^{\frac{m}{4}} & \mbox{si} & \mbox{$m/2$ est pair} \\
\PPhi(\tilde{x}_{1-\frac{m}{2}}) & = & \om^{\frac{m}{2}}
(\tp^{-1} \sp^{-1})^{\frac{m}{4}} (\tp \sp)^{\frac{m}{4}} & \mbox{si} & \mbox{$m/2$ est pair} \\
\PPhi(\tilde{x}_{\frac{m}{2}}) & = & \om^{\frac{m}{2}}
(\underbrace{\sp^{-1} \tp^{-1}\dots }_{\frac{m}{2}})
(\underbrace{\tp \sp \dots }_{\frac{m}{2}}) & \mbox{si} & \mbox{$m/2$ est impair} \\
\PPhi(\tilde{x}_{1-\frac{m}{2}}) & = & \om^{\frac{m}{2}}
(\underbrace{\tp^{-1} \sp^{-1}\dots }_{\frac{m}{2}})
(\underbrace{\sp \tp \dots }_{\frac{m}{2}}) & \mbox{si} & \mbox{$m/2$ est impair} \\
\end{array} 
$$
\caption{Valeurs de $\PPhi$ sur la famille $\tilde{x}$}
\label{tablephix}
\end{table}
Soit enfin $\Phi' =\J(\Phi^{-1})$. Comme $J(\tilde{x}_r) = \tilde{x}_{1-r}$
et $\Ad(\Phi) \circ \J$ envoie
$\sigma_{\Phi'}$ et $\tau_{\Phi'}$ respectivement sur $\tp$
et $\sp$, on a
\begin{eqnarray}
\label{phio} 
\Ad(\Phi) \circ \J \circ \PPhi' = \PPhi \circ J
\end{eqnarray}
\begin{figure}
\begin{center}
\includegraphics{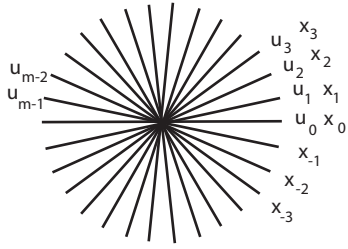}
\end{center}
\caption{$(u_r)$ et $(x_r)$}
\end{figure}

\begin{table}
$$
\begin{array}{|l|l|}
\hline
\multicolumn{2}{|l|}{\mbox{Pour tout $m$ : }} \\
\hline
\PPhi(\tilde{u}_{2r})    =  (\tp \sp)^r \bullet \sp^2 
& 2 \leq 2r < \frac{m}{2} \\
\hline
\PPhi(\tilde{u}_{2r+1})    =  (\tp \sp)^r \bullet \tp^2 
& 2 \leq 2r+1 < \frac{m}{2} \\
\hline
\multicolumn{2}{|l|}{\mbox{Pour $m$ impair : }} \\
\hline
\PPhi(\tilde{u}_{2r})  =  (\tp \sp )^{\frac{2r-1-m}{2}} \bullet \tp^2  & \frac{m+1}{2} < 2r \leq m-1 \\
\hline
\PPhi(\tilde{u}_{2r+1})  =  (\tp \sp )^{\frac{2r+1-m}{2}} \bullet \sp^2  & \frac{m+1}{2} < 2r+1 \leq m-1 \\
\hline
\multicolumn{2}{|l|}{\mbox{Pour $m$ pair : }} \\
\hline
\PPhi(\tilde{u}_{m-2r})  =  (\tp \sp)^{-r} \bullet \sp^2 &
 1 \leq 2r < \frac{m-2}{2} \\
\hline
\PPhi(\tilde{u}_{m-2r+1})  =  (\tp \sp)^{-r} \bullet \tp^2 &
 1 \leq 2r+1 < \frac{m-2}{2}  \\
\hline
\PPhi(\tilde{u}_{\frac{m}{2}})  =  \om^{\frac{m}{2}}
(\sp^{-1} \tp^{-1})^{\frac{m}{4}} \sp^2 (\sp^{-1} \tp^{-1})^{\frac{m}{4}}
 & \mbox{$\frac{m}{2}$ pair} \\
\hline
\PPhi(\tilde{u}_{\frac{m}{2}})  =  \om^{\frac{m}{2}}
(\underbrace{\sp^{-1} \tp^{-1} \dots }_{\frac{m}{2}})
\tp^2
(\underbrace{\tp^{-1} \sp^{-1} \dots }_{\frac{m}{2}})
 & \mbox{$\frac{m}{2}$ impair} \\
\hline
\PPhi(\tilde{u}_{\frac{m}{2}+1})  =  \om^{\frac{m}{2}}
(\sp^{-1} \tp^{-1})^{\frac{m}{4}} \tp^2 (\tp^{-1} \sp^{-1})^{\frac{m}{2}}
(\tp \sp)^{\frac{m}{4}}  & \mbox{$\frac{m}{2}$ pair} \\
\hline
\PPhi(\tilde{u}_{\frac{m}{2}+1})  =  \om^{\frac{m}{2}}
(\sp^{-1} \tp^{-1})^{\frac{m-2}{4}} (\tp^{-1} \sp^{-1})^{\frac{m}{2}}
\sp^2 (\tp \sp)^{\frac{m-2}{4}}
 & \mbox{$\frac{m}{2}$ impair} \\
\hline
\end{array}
$$
\caption{Valeurs de $\PPhi$ sur la famille $\tilde{u}$}
\label{tablephiu}
\end{table}

\subsection{Existence et unicit\'e d'une solution}

Pour tout $n \geq 0$ on note $\mathcal{A}_n'(\k)$
l'id\'eal de $\mathcal{A}'(\k)$ form\'e des \'el\'ements de valuation
au moins $n$.

Notons $\tilde{u}_i = \exp v_i$, $v_i \in L'(\k)$, et $\dd \PPhi :
L'(\k) \to \mathcal{A}'(\k)$ l'application tangente.
On d\'eduit des formules
de la table \ref{tablephiu} que, pour tout $0 \leq i \leq m-1$, on a $\dd \PPhi(v_i) \equiv 2
t_i \mod \mathcal{A}'_2(\k)$. C'est immédiat pour
$m$ impair (pour $i=m-1$ on utilise
$\tilde{u}_{m-1} \dots \tilde{u}_1 \tilde{u}_0 = 1$), et pour $m$
pair les seuls cas non triviaux sont $i \in \{ \frac{m}{2}, 1 + \frac{m}{2} \}$.
Dans ces cas, on remarque
$(\tp \sp)^{\frac{m}{2}} \equiv \om^{\frac{m}{2}}$ modulo
$\mathcal{A}'_2(\k)$, qui d\'ecoule de ce que $\sum_{k=1}^{\frac{m}{2}}
\om^{-k} - \sum_{k=0}^{\frac{m}{2}-1} s \om^k \in \k W$ agit
par $0$ sur la composante homogène de degré $1$ de $\mathcal{A}'(\k)$.

En particulier, si $\Phi^1 \equiv \Phi^2 \mod
\mathcal{A}'_r(\k)$ on a $\dd \PPhi^1(v_i) \equiv \dd \PPhi^2(v_i)
\mod \mathcal{A}'_{r+1}(\k)$. On peut alors consid\'erer
l'\'equation (\ref{maineq}) pour $\Phi \in \mathcal{A}'(\k)$
tel que $\Phi \equiv 1 \mod \mathcal{A}'_1(\k)$.

\begin{prop} \label{propexsol} Il existe un
unique $\Phi \equiv 1 \mod
\mathcal{A}_1'(\k)$ dans $\mathcal{A}'(\k)$ qui v\'erifie
(\ref{maineq}). De plus, un tel $\Phi$ appartient \`a
$\exp \mathfrak{g}'(\k)$.
\end{prop}
\begin{proof}
On a $\partial \Phi,\dd \PPhi(\psi) \in \mathcal{A}_1'(\k)$
donc (\ref{maineq}) est v\'erifi\'ee modulo $\mathcal{A}_1(\k)$.
Supposons d\'emontr\'e qu'existe $\Phi \equiv 1 \mod
\mathcal{A}'_1(\k)$ v\'erifiant (\ref{maineq}) modulo
$\mathcal{A}'_r(\k)$ pour un certain $r \geq 1$, et que toute autre
solution soit congrue \`a $\Phi$ modulo $\mathcal{A}'_r(\k)$. Pour
montrer la m\^eme chose au rang $r+1$, il suffit de montrer
qu'existe un unique $\Phi' \in \mathcal{A}'(\k)$ homog\`ene
de degr\'e $r$ tel que $\Phi + \Phi'$ soit solution de (\ref{maineq})
modulo $\mathcal{A}'_{r+1}(\k)$, c'est-\`a-dire
$$
 (\partial \Phi + \partial \Phi') (\Phi + \Phi')^{-1}\equiv
- \dd \left( \widetilde{ \Phi + \Phi'} \right) (\psi) \mod
\mathcal{A}'_{r+1}(\k).
$$
Or $\dd \widetilde{\Phi + \Phi'} (v_i) \equiv \dd \PPhi(v_i) \mod
\mathcal{A}'_{r+1}(\k)$ donc $\dd \widetilde{\Phi + \Phi'}
(\psi) \equiv \dd \PPhi(\psi) \mod \mathcal{A}'_{r+1}(\k)$. D'autre
part
$\partial \Phi + \partial \Phi' \in \mathcal{A}'_1(\k)$
donc
$$
 (\partial \Phi + \partial \Phi') (\Phi + \Phi')^{-1} \equiv
 (\partial \Phi + \partial \Phi') \Phi^{-1} \mod
\mathcal{A}'_{r+1}(\k),
$$
et l'\'equation est \'equivalente \`a $(\partial \Phi
+ \partial \Phi' )\Phi^{-1} \equiv - \dd \PPhi(\psi) \mod \mathcal{A}'_{r+1}(\k)$
c'est-\`a-dire $\partial \Phi' \equiv - \left( \dd \PPhi(\psi) \right)
\Phi - \partial \Phi$ et $\partial \Phi'$ est uniquement
d\'etermin\'e. Enfin, $\partial$ se restreint en une bijection sur
chacune des composantes homog\`enes de degr\'e au moins 1 de
$\mathcal{A}'(\k)$, donc $\Phi'$ existe et est uniquement
d\'etermin\'e, ce qui nous permet de conclure par r\'ecurrence
l'existence et l'unicit\'e d'une solution.

Comme $\psi \in \mathcal{G}(\k)$ et d'apr\`es les
formules de la table \ref{tablephiu} on a $\dd \PPhi (\psi) \in \exp 
\mathfrak{g}'(\k)$ donc $\Delta(\dd \PPhi(\psi)) = \dd \PPhi(\psi)
\otimes \dd \PPhi(\psi)$. 
Par l'action naturelle de $\kt$ sur $\mathcal{A}'(\k) \widehat{\otimes}
\mathcal{A}'(\k)$ on peut définir, de façon analogue à $\partial$,
une dérivation de $\mathcal{A}'(\k) \widehat{\otimes}
\mathcal{A}'(\k)$ que l'on note encore $\partial$. On a 
$$
\partial ( \Phi \otimes \Phi) = \partial \Phi \otimes \Phi
+ \Phi \otimes \partial \Phi = - ( \dd \PPhi(\psi) \otimes 1
+ 1 \otimes \dd \PPhi(\psi))(\Phi \otimes \Phi)
$$
soit $\partial(\Phi \otimes \Phi) = - \Delta(\dd \PPhi(\psi))(\Phi \otimes
\Phi)$.

Comme $\Delta \circ \partial$ et $\partial \circ \Delta$
sont deux applications lin\'eraires continues
$g:\mathcal{A}'(\k) \to \mathcal{A}'(\k) \widehat{\otimes}
\mathcal{A}'(\k)$ qui v\'erifient
toutes deux l'\'equation $g(RS) = g(R) \Delta(S) + \Delta(R) g(S)$
et coïncident sur les g\'en\'erateurs $(t_i)$, on a
$\Delta \circ \partial = \partial \circ \Delta$ et 
$\partial(\Delta(\Phi)) = \Delta(\partial \Phi) = - \Delta(
\dd \PPhi(\psi)) \Delta(\Phi)$. Ainsi $\Phi \otimes \Phi$ et
$\Delta(\Phi)$ ont la m\^eme image par l'application
$K \mapsto (\partial K) K^{-1}$, donc ils sont \'egaux
parce qu'ils ont m\^eme terme constant, donc $\Phi$ est grouplike
c'est-\`a-dire $\Phi \in \exp \mathfrak{g}'(\k)$.

\end{proof}

\begin{ccor} Si $\Phi \equiv 1 \mod
\mathcal{A}'_1(\k)$ et $\Phi$ v\'erifie (\ref{maineq}), alors
$\Phi^{-1} = \J( \Phi)$.
\end{ccor}
\begin{proof}
On pose $\Phi' = \J(\Phi^{-1})$. D'apr\`es (\ref{phio}) on a $\Ad(\Phi) \circ \J \circ \dd \PPhi'
= \dd \PPhi \circ J$. D'autre part, comme $\psi \in \mathcal{G}'(\k]$
on a $J(\psi)=-\psi$
Alors $\Phi'\J(\Phi)  = 1$
implique
$$
D(\Phi') = - \Phi' \bullet \J(D(\Phi)) = - \Phi' \bullet
\J(\dd \PPhi ( - \psi )) = - \Phi' \bullet \J \circ \dd \PPhi \circ J(\psi)
$$
Or $\Ad(\Phi') \circ \J \circ \dd \PPhi \circ J = \Ad(\Phi') \circ \J \circ
\Ad(\Phi) \circ \J \circ \dd \PPhi' = 
\Ad(\Phi' \J(\Phi)) \circ \dd \PPhi' = \dd \PPhi'$
d'o\`u $\Phi'$ v\'erifie (\ref{maineq}) d'o\`u $\Phi = \Phi'$ et
$\Phi^{-1} = \J (\Phi)$ par unicit\'e.

\end{proof}

Il reste à vérifier qu'une telle solution vérifie l'équation du
demi-tour.

\subsection{L'alg\`ebre de Lie $\mathcal{G}'(\k)$}
\label{ssalglie}
Soit $(l,\psi) \in \mathcal{G}'(\k)$. Comme $J(f) = f^{-1}$ pour tout
$(\la,f) \in G'(\k)$ on a, d'une part $J(\psi)  = -\psi$, d'autre part
les \'equations (\ref{gexp2}) se simplifient. En effet, posons
$F = \sss f^{-1}$. On a $L = F ((\tau \sigma)^{\frac{m+1}{2}} \bullet F)
= F (\tau \bullet J(F))$. De même $M = (\tau^{-1} \bullet F) J(F)$
d'où $(\sigma \tau)^{-k} \bullet M = \tau^{-1}(\tau \sigma)^{-k}
\bullet (F (\tau \bullet J(F)))$. Posons $F' = F (\tau \bullet J(F))$.
Puisque $\tau$ commute à $\ooo$
lorsque $m$ est pair, on en d\'eduit que $f$ vérifie les équations

\begin{eqnarray}
\label{gexp2b}
\left( \prod_{k= \frac{m-1}{2}}^1 (\tau \sigma)^{-k} \bullet F')
 \right) F = \ooo, \prod_{k = \frac{m}{2} -1 }^{0} (\tau \sigma)^{-k}
\bullet F' = \ooo
\end{eqnarray}
suivant que $m$ est pair ou impair. 

Pour lin\'eariser (\ref{gexp2b}) on pose $\la = 1+\eps l$,
$ f= \exp(\eps \psi)$ avec $\psi \in L'(\k)$, $\eps^2 =0$. On a
$\sigma^2 = u_0 = \exp(v_0)$ d'o\`u $\sss = \exp (\eps \frac{l}{2} v_0)$,
$\ooo = \exp( \eps \frac{l}{2} z)$. On en d\'eduit
$
F = \sss f^{-1} = \exp \eps (\frac{l}{2} v_0 -\psi)$
et $\tau \bullet J(F) = \exp \eps(\frac{l}{2} v_1 - \tau\bullet J(\psi))$.
Soient $X_0,\dots,X_{m-1}$ une famille d'indéterminées indexée par $\Z/m\Z$,
et $\mathcal{R}= \k\ll X_0,\dots,X_{m-1}\gg/(e^{X_{m-1}}\dots e^{X_0}-1)$.
On note 
$\eta$ l'automorphisme de $\mathcal{R}$ défini par $X_i \mapsto
X_{i+1}$, et pour $S \in \mathcal{R}$ on note
$\eta.S$ l'image de $S$ par cet endomorphisme.
Pour tout $m$-uplet
$\underline{x} = (x_0,\dots,x_{m-1})$ d'éléments de
$L'(\k)$ ou $\mathcal{A}'_1(\k)$ tels que $e^{x_{m-1}}\dots e^{x_0}=1$
on définit $S(\underline{x})$ par
spécialisation. Soit $\underline{\tilde{v}} = (\tilde{v}_0,
\dots, \tilde{v}_{m-1})$, avec $\tilde{v}_i = v_i - \frac{z}{m}$ et $S \in \mathcal{R}$ telle que
$\psi = S(\underline{\tilde{v}})$. 

On a $(\tau \sigma)\bullet \tilde{v}_r = \tilde{v}_{r+2}$ et on
vérifie facilement $\tau^{-1} \bullet u_r = J(u_{r-1})$
donc $\tilde{v}_r = \tau \bullet J(\tilde{v}_{r-1})$. On en déduit
que $\tau \bullet J(\psi) = (\eta.S)(\underline{\tilde{v}})$
et $\tau \sigma \bullet \psi = (\eta^2.S)(\underline{\tilde{v}})$.
On déduit alors de (\ref{gexp2b}) que
\begin{eqnarray}
\label{lieeq}
 \sum_{k=0}^{m-1} \eta^k.S(\underline{\tilde{v}}) - \frac{l}{2} 
\sum_{k=0}^{m-1} \tilde{v}_k = 0
\end{eqnarray}
On en déduit que, pour tout $m$-uplet $\underline{d} \in
\mathcal{A}'_1(\k)^m$ tel que $e^{d_{m-1}}\dots e^{d_0}\equiv 1 \mod
\mathcal{A}_{n+1}(\k)$,
\begin{eqnarray}
\label{liegimpmod}
 \sum_{k=0}^{m-1} \gamma^k.S(\underline{d}) - \frac{l}{2} 
\sum_{k=0}^{m-1} d_k \equiv 0 \mod \mathcal{A}'_{n+1}(\k)
\end{eqnarray}

\subsection{L'\'equation du demi-tour}

On suppose que $\Phi$ satisfait (\ref{maineq}), et on pose
$\PP = \Phi e^{t_1}$, $\xi = \Phi e^{t_1} (\om s \bullet \Phi^{-1}) e^{t_2}$.
Les équations (\ref{impeq3}) et (\ref{paireq3}) s'écrivent alors
$Q = 1$ avec
$$
Q = \prod_{r=0}^{\frac{m}{2}-1} (\om^r \bullet \xi),
Q = \left[\prod_{r=0}^{\frac{m-3}{2}} (\om^r \bullet \xi)\right] (\om^{\frac{m-1}{2}} \bullet \PP)
$$
suivant que $m$ est pair ou impair. On a $\PP \om s \Phi = \tp$
et $\tp \sp = \xi \om$ d'où 
$$
\prod_{r=0}^{n-1} (\om^r \bullet \xi) = (\tp \sp)^n \om^{-n}.
$$
Posons $D(Q) = (\partial Q) Q^{-1}$. De $D(AB) = D(A) + A \bullet D(B)$
on déduit $D(\Phi^{-1}) = - \Phi \bullet D(\Phi)$
et $D(\xi) = D(\Phi) + \PP \bullet t_1 - \tp \bullet D(\Phi)
+ \xi \bullet t_2$. On note $q_0 = 1$, $q_1 = \PP$, $q_{2r} = (\tp \sp)^r$,
$q_{2r+1} =  (\tp \sp)^r \PP$, et $t'_r = q_r \bullet t_r$. Ainsi
$$
D(Q) = \sum_{r=1}^{m} t'_r + \sum_{r=0}^{\frac{m-1}{2}}
(\tp \sp)^r \bullet D(\Phi) - \sum_{r=0}^{\frac{m}{2}-1}
(\tp \sp)^r \tp \bullet D(\Phi)
$$
que $m$ soit pair ou impair. On veut montrer $Q = 1$. Pour ce faire,
on va montrer par récurrence $Q \equiv 1 \mod \mathcal{A}'_n(\k)$
pour tout $n \geq 1$. On note que, à $n$ fixé, c'est équivalent à
$D(Q) \equiv 0$. Soit donc $n \geq 1$. A partir de maintenant, les
congruences seront toujours comprises modulo $\mathcal{A}_{n+1}(\k)$.
On suppose donc que $Q \equiv 1 +R$ avec $R$ polynôme de Lie
homogène de degré $n$.
En utilisant les formules de la table (\ref{tablephiu}) on en déduit par
un calcul direct les formules
\begin{eqnarray}
\forall r\geq 0 \ \ t'_{r+m} \equiv t'_r \\ 
\forall r\in [0,m-1]\ \ \dd \PPhi(\tilde{v}_r) \equiv 2 t'_r. 
\end{eqnarray}
D'autre part, en utilisant les identités
$q_{2r+1}^{-1} q_{2r} = \om^2 \bullet (e^{-t_1} \Phi^{-1})$
et $q_{2r}^{-1} q_{2r-1} = \om^r s \bullet (e^{-t_0} \Phi)$,
on montre par récurrence sur $r \in [1,m-1]$ la formule 
\begin{eqnarray}
e^{2t'_{m-1}} e^{2t'_{m-2}} \dots e^{2t'_{m-r}}
\equiv Q e^{-t_0} (s \bullet q_r ) e^{t_{m-r}} q_{m-r}^{-1}
\end{eqnarray}
En particulier, pour $r = m-1$ on a
$$
e^{2 t'_{m-1}} \dots e^{2 t'_1} \equiv Q e^{-t_0} (s \bullet Q) e^{-t_0}
$$
On pose $t''_{m-1} = t'_{m-1} - R$, $t''_1 = t'_1 - (s \bullet R)$,
$t''_i = t'_i$ si $i \not\in \{1, m-1\}$ et $\underline{t}'' =
(t''_0,\dots,t''_{m-1})$, de telle façon que
$e^{2 t''_{m-1}} \dots e^{2 t''_0} \equiv 1$. On utilise maintenant
$$
D(\Phi) = - \dd \PPhi(\psi) = - \dd \PPhi( S( \underline{v})) = - S(\dd \PPhi(
\underline{v})) \equiv - S(
\underline{t}'')
$$
On en d\'eduit
$(\tp \sp)^k \bullet D(\Phi) \equiv - S( (\tp \sp)^k \bullet \underline{t}'')
\equiv - (\eta^{2k} .S)(\underline{t}'')$.
D'autre part,
$$
\begin{array}{lclcl}
 \tp \bullet D(\Phi) & = & - \tp \bullet \dd \PPhi(\psi)
& = & \tp \bullet \dd \PPhi \circ J(\psi) \\
& = & \Ad(\tp \Phi) \circ \J \circ \dd \PPhi(\psi)
& = & \Phi \om s e^{t_1} \bullet \J(\dd \PPhi(\psi)) \\
& = & \J ( \Phi^{-1} s e^{t_0} \bullet \dd \PPhi(\psi)) \\
\end{array}
$$
Comme on d\'eduit des d\'efinitions que $\Phi^{-1} s e^{t_0} \bullet
t'_r = \J (t'_{r+1})$, on a
$$
 \tp \bullet D(\Phi) \equiv \J ( \Phi^{-1} s e^{t_0} \bullet
S( \underline{t}'')) \equiv (\eta .S)(\underline{t}'')
$$
et plus g\'en\'eralement $ (\tp \sp)^k \tp \bullet D(\Phi) \equiv
(\eta^{2k+1} .S)(\underline{t}'')$.
On a donc
$$
D(Q) \equiv t'_1 + \dots + t'_m - \sum_{k=0}^{m -1}
(\eta^k . S)(\underline{t}'')
$$
Appliquant (\ref{liegimpmod}) avec $l = 1$ on en déduit
$D(Q) \equiv R + s \bullet R$. En utilisant $\J(\Phi) = \Phi^{-1}$ si $m$
est impair et la centralité de $\om^{\frac{m}{2}}$ dans
le cas pair, on a $\om \bullet Q =  
\xi^{-1} \bullet Q$ donc $\om \bullet R = R$. Ainsi la classe de $D(Q)$, donc $R$, est invariant
sous l'action de $W$, et 
$D(Q) \equiv n R \equiv 2R$.
Ainsi $R =0 $, sauf peut-\^etre si $n =2$.
Pour conclure il nous suffirait
de montrer que l'espaces des polyn\^omes de Lie homogènes de degré 2
n'admet pas
de vecteur invariant pour l'action de $W$. C'est vrai si $m$
est impair.

En effet, notons $V$ (resp. $V'$) l'espace des \'el\'ements homog\`enes de degr\'e
1 de $\mathcal{A}(\k)$ (resp. $\mathcal{A}'(\k)$) muni de l'action
restreinte de $W$. Sous cette action,
$V \simeq V' \oplus \un$, et $\Lambda^2 V \simeq \Lambda^2 V' \oplus V'$. Il
suffit de montrer que $E = \Lambda^2 V$ ne contient pas la
repr\'esentation triviale de $W$. On constate ais\'ement que
le caract\`ere $\chi_E$ de $E$ vaut $0$ sur les rotations non
triviales et $\frac{1-m}{2}$ sur les r\'eflections. Comme $\dim E
= \frac{m(m-1)}{2}$ on en d\'eduit que le produit scalaire de $\chi_E$
avec $\un$ est nul, ce qui conclut.

Dans le cas pair, il nous faut utiliser de plus l'invariance par $\J$.
Comme $\J (\Phi) = \Phi^{-1}$ on a $\J (\xi ) = \Phi^{-1} s e^{t_0}
\bullet \xi$ et $\J(Q) = (\Phi^{-1} e^{t_0}) \bullet (s \bullet Q)$
donc $\J (Q) \equiv s \bullet Q \equiv Q$ et $R$ est invariant sous l'action
de $\J$. La conclusion d\'ecoule alors du
fait que
l'espace des polyn\^omes de Lie homog\`enes de degr\'e 2 de $\mathcal{A}'(\k)$
invariants par $W$ et $\J$ est nul.

En effet, 
si l'on note $\Omega = \Ad(s) \circ \J \circ \Ad(\om)
$, $S = \Ad(s)$. Alors $\Omega(t_r) = t_{r+1}$, $\Omega \circ S \circ \Omega
= S$, $S^2= 1$. Soit $W'$ le groupe di\'edral engendr\'e par $\Omega$
et $S$. Il suffit de montrer que $E = \Lambda^2 V$ ne contient pas
de vecteur invariant par $W'$. Le caract\`ere $\chi_E$ de $E$ vaut 0 sur
les rotations non triviales, $\chi_E(\Omega^{\frac{m}{2}}) = -m/2$, 
$\chi_E(S) = (4-m)/2$ et $\chi_E(S \Omega^{\frac{m}{2}}) = (2-m)/2$.
On en d\'eduit $(\chi_E \mid \un) = 0$, ce qui conclut.

\section{Appendice 1 : \'equations Fuchsiennes formelles}
\setcounter{equation}{0}
Nous établissons ici des résultats bien connus sur les équations
différentielles fuchsiennes à valeurs dans une
algèbre de Hopf complète, pour lesquels nous n'avons pas trouvé de
référence convenable.

Soit $\IXE$ un ensemble fini, et $M$ le mono\"\i de libre sur $\IXE$.
On note $\emptyset$ son \'el\'ement neutre, et $(\alpha,\beta)
\mapsto \alpha \star \beta$ son produit de concatenation.
Pour tout corps $\k$ on
note $\mathcal{M}(\k)$ la $\k$-alg\`ebre de Magnus sur $\IXE$,
c'est-\`a-dire l'ensemble des s\'eries formelles sur $\k$ en les
variables non commutatives \'el\'ements de $\IXE$. Si $\k$ est
un corps topologique, on munit $\mathcal{M}(\k)$ de la
topologie produit, c'est-\`a-dire la topologie de la convergence
simple, et de sa graduation naturelle $\mathcal{M}(\k)
= \prod \mathcal{M}_{\alpha}(\k)$ et des projections
naturelles $\pi_{\alpha} : \mathcal{M}(\k) \to \mathcal{M}_{\alpha}(\k)= \k$
pour $\alpha \in M$. Enfin, on note $\omega : \mathcal{M}(\k) \to \N$
la valuation associ\'ee \`a la graduation totale ($\omega(x) = 1$
pour $x \in \IXE$).

Pour tout ouvert $U \subset \C$ (resp. $U \subset \R$) on dit classiquement
que
$C : U \to \mathcal{M}(\C)$ (resp. $C : U \to \mathcal{M}(\R)$)
est analytique sur $U$ si et seulement si $\pi_{\alpha} \circ C$
est analytique sur $U$ pour tout $\alpha \in M$. 

Soient $A \in \IXE$, $D = \{ z \in \C \ \mid \ |z| < 1 \}$ et
$C : D \to \mathcal{M}(\C)$ analytique sur $D$ telle que,
pour tout $z \in D$, $C(\bar{z}) = \overline{C(z)}$ et
$\omega(C(z)) \geq 1$. On consid\`ere l'\'equation diff\'erentielle
\begin{eqnarray}
G'(z) & = & \left( \frac{A}{z} + C(z) \right) G(z) \label{eqgencomp}
\end{eqnarray}
Elle se restreint en une \'equation diff\'erentielle r\'eelle
sur $]0,1[$,
\begin{eqnarray}
G'(x) & = & \left( \frac{A}{x} + C(x) \right) G(x) \label{eqgenreel}
\end{eqnarray}

\begin{lemme} \label{unicitepositive} Il existe une seule solution $G_+$
de (\ref{eqgenreel}) sur $ ]0,1[$ telle que $G_+(x) \sim x^A$ quand $x \to 0^+$.
\end{lemme}
Dans cet \'enonc\'e, $x^A = \exp(A \log(x))$, et $G(x) \sim x^A$ signifie
que $G(x) x^{-A}$ est une fonction
analytique sur $]-1,1[$, qui vaut 1 en 0.
\begin{proof}
On cherche $G_+(x)$ sous la forme $P(x) x^A$ avec $P$ analytique
et $P(0) = 1$. Ainsi $P(x) = \sum_{n=0}^{\infty}
p_n x^n$, $p_0 = 0$ et $p_n \in \mathcal{M}(\R)$.
L'\'equation (\ref{eqgenreel}) se r\'e\'ecrit sous la forme
$$
P(x) \frac{A}{x} + P'(x) = \frac{A}{x} P(x) + C(x) P(x)
$$
soit $(*) \ \ x P'(x) + [P(x),A] = x C(x) P(x)$. Notons
$C(x) = \sum_{n=0}^{\infty} c_n x^n$. L'\'equation implique alors
$$
n p_n+ [p_n,A] = \sum_{r+s = n} c_{r+1} p_s
$$
pour tout $n \geq 0$. Le membre de droite ne comporte que des $p_s$
pour $s < n$, donc l'unicit\'e et l'existence d'une solution formelle
d\'ecoulent par r\'ecurrence de l'inversibilit\'e dans
$\mathcal{M}(\k)$, pour tout $n \geq 1$, de l'op\'erateur
$n - \mathrm{ad}(A)$. Il reste \`a v\'erifier l'analyticit\'e de $P$
ainsi d\'efini.

Notons
$$
P(x) = \sum_{\alpha \in M} f_{\alpha}(x) \alpha, \ \ \ 
C(x) = \sum_{\alpha \in M} g_{\alpha}(x) \alpha \ \ \ 
$$
avec $f_{\alpha}, g_{\alpha} \in \R[[x]]$, les $g_{\alpha}$
\'etant analytiques. Comme $P(0) = 0$, on a $f_{\alpha}(0) = 0$
si $\omega(\alpha) \geq 1$, $f_{\emptyset}(0) = 1$. L'\'equation
$(*)$ s'\'ecrit alors
$$
\sum_{\alpha \in M} \left( x f_{\alpha}'(x) \alpha + f_{\alpha}(x)(\alpha \star A
- A \star \alpha) - \sum_{\beta \star \gamma = \alpha} x g_{\beta}(x)
f_{\gamma}(x) \alpha \right) = 0.
$$
En particulier $x f_{\emptyset}'(x) = x g_{\emptyset}(x) f_{\emptyset}(x)$
donc $f_{\emptyset}'(x) = g_{\emptyset}(x) f_{\emptyset}(x)$
pour $x \neq 0$. Or $g_{\emptyset}(x) = 0$ donc $f_{\emptyset}'(x) =0$
et comme $f_{\emptyset}(0) = 1$, on a $f_{\emptyset}(x) = 1$.

Pour tout $\alpha \in M$, on note $\alpha'$ (resp. $\alpha''$) l'unique
mon\^ome tel que $\alpha = \alpha' \star A$ (resp. $\alpha = A \star
\alpha''$) s'il existe, $\dagger$ sinon et on pose $f_{\dagger}(x) = 0$.
Comme $g_{\emptyset}(x) = 0$ on a, pour tout $\alpha$ tel que $\omega
(\alpha) \geq 1$,
$$
f_{\alpha}'(x) = \frac{f_{\alpha''}(x) - f_{\alpha'}(x)}{x}
+ \sum_{\stackrel{\beta \star \gamma = \alpha}{\omega(\gamma)< \omega
(\alpha)}} g_{\beta}(x) f_{\gamma}(x).
$$
Supposons que les $f_{\gamma}(x)$ pour $\omega(\gamma) < \omega(\alpha)$
soient analytiques sur $]-1,1[$. De plus $f_{\gamma}(0)
= 0$, sauf si $\gamma = \emptyset$. Ainsi $f_{\alpha''}(0) = 
f_{\alpha'}(0) = 0$ sauf si $\alpha = A$. Mais alors $\alpha' = \alpha''
= \emptyset$ et  $f_{\alpha''}(x) - f_{\alpha'}(x) = 0$. Dans tous
les cas on en d\'eduit que $f_{\alpha}'(x) \in \R[[x]]$ et admet
un rayon de convergence au moins \'egal \`a 1.
Il en est donc de m\^eme pour $f_{\alpha}$,
et on conclut par r\'ecurrence sur $\omega(\alpha)$.
\end{proof}

On d\'emontre de la m\^eme fa\c con le lemme 
\begin{lemme} \label{unicitenegative} Il existe une seule solution $G_-$
de (\ref{eqgenreel}) sur $]-1,0[$ telle que $G_-(x) \sim (-x)^A$ quand $x \to 0^-$.
\end{lemme}
\begin{figure}
\begin{center}
\includegraphics{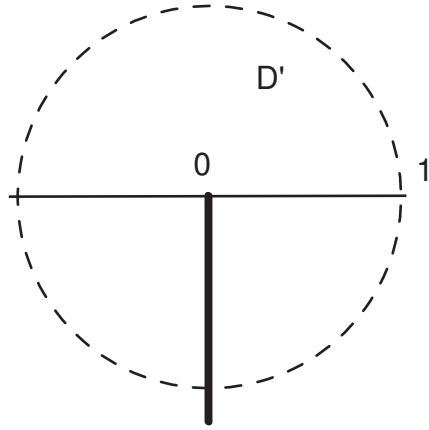}
\end{center}
\caption{Le domaine $D'$} \label{genmono}
\end{figure}
On montre facilement que ces solutions $G_+$, $G_-$ s'\'etendent
analytiquement au domaine complexe simplement connexe
$D'= D \setminus \ii \R_-$ (voir figure \ref{genmono}) en des solutions de l'\'equation (\ref{eqgencomp}).
Pour \'etudier ces prolongements, on introduit $\Log(z)$, branche
du logarithme complexe d\'efini comme le prolongement de $\log(x)$
pour $ x\in \R_+^*$ sur $D'$ et, pour
$z \in D'$, on note $z^A = \exp(A \Log(z))$.
On a,
pour $x \in ]0,1[$, $\Log(-x) = \log(x) + \ii \pi$ d'o\`u
$$
(-x)^A = \exp(A \Log(-x)) = \exp\left( A (\ii \pi + \log x) \right) =
x^A \exp(\ii \pi A).
$$
Notons alors $\tilde{G}(z) = G_+(z) e^{-\ii \pi A}$. C'est une fonction
analytique sur $D'$ qui v\'erifie l'\'equation (\ref{eqgencomp}).
Il existe une unique fonction $P_+(x)$ analytique r\'eelle sur $]-1,1[$
telle que $G_+(x) = P_+(x) x^A$ pour tout $x \in ]0,1[$. On en
d\'eduit $G_+(z) = P_+(z) z^A$ pour tout $z \in D'$. En
particulier, pour $x \in ]-1,0[$,
$$
\tilde{G}(x) = P_+(x) x^A e^{-\ii \pi A} = P_+(x) (-x)^A. 
$$
On en d\'eduit que $\tilde{G}(x)$ est une fonction
analytique r\'eelle sur $]-1,0[$ telle que $\tilde{G}(x) \sim (-x)^A$ quand
$x$ tend vers $0^-$. D'apr\`es le lemme \ref{unicitenegative} cela
signifie
\begin{lemme} Pour tout $x \in ]-1,0[$, $G_+(x) = G_-(x) e^{\ii \pi A}$
\label{passage}
\end{lemme}

L'algèbre $\mathcal{M}(\k)$ est la complétion de l'algèbre enveloppante
de l'algèbre de Lie libre sur $\IXE$. Elle est donc munie
d'une structure d'algèbre de Hopf complète au sens de \cite{QUILLEN}.
Soit $V \subset \Q \IXE \subset \mathcal{M}(\k)$ un sous-espace vectoriel
des combinaisons linéaires rationnelles d'éléments de $\IXE$.
On note $\mathcal{U}(\k)$ le quotient de $\mathcal{M}(\k)$ par
l'idéal de Hopf fermé qu'il engendre, et $\mathcal{U}^1(\k) = 
\k \IXE / V \otimes \k \subset \mathcal{U}(\k)$. Soient $A$
l'image dans $\mathcal{U}(\k)$ d'un $\hat{A} \in \IXE$ et
$C : D \to \mathcal{U}^1(\C)$ tel que $C(\bar{z}) = \overline{C(z)}$.
Si l'on considère l'équation \ref{eqgenreel} avec cette fois $G$
à valeurs dans $\mathcal{U}(\R)$, les lemmes précédents sont encore
valables. En effet, l'unicité des solutions se démontre de la même
façon, et leur existence découle d'un relèvement de $A$ en
$\hat{A}$ et de $C$ en une fonction analytique $\hat{C} : D \to
\mathcal{M}(\C)$. On invoquera donc ces lemmes également
dans ce cadre.

Soit $\Delta$ le coproduit de $\mathcal{U}(\k)$, et notons
$B(z) = \frac{A}{z} + C(z)$.
Pour tout $z \in D'$, comme $C(z)$ est primitif $B(z)$ l'est également.
Alors 
$$
\Delta(G(z))' = \Delta(G'(z)) = \Delta(B(z)) \Delta(G(z)) = (B(z) \otimes 1 + 1 \otimes
B(z))\Delta(G(z))
$$
Soit $\hat{G}(z) = G(z) \otimes G(z)$. On a
$$
\hat{G}'(z) = G'(z) \otimes G(z) + G'(z) \otimes G(z) = (B(z) \otimes 1
+ 1 \otimes B(z)) \hat{G}(z)
$$
Ainsi, $G(z) \otimes G(z)$ et $\Delta(G(z))$ v\'erifient
la m\^eme \'equation diff\'erentielle, et la m\^eme condition
asymptotique parce que $z^A$ est grouplike. Le lemme \ref{unicitepositive}
appliqué à $\mathcal{U}(\k) \widehat{\otimes} \mathcal{U}(\k)$
considéré comme quotient de l'algèbre de Magnus sur $\IXE \times \IXE$
permet de conclure la démonstration du lemme suivant :

\begin{lemme} Si $C(x)$ est primitif pour tout $x \in ]-1,1[$
alors, pour tout $z \in D'$, $G_+(z)$ et $G_-(z)$ sont
grouplike. \label{grouplike}
\end{lemme}

\section{Appendice 2 : le cas $m =3$ et les associateurs de Drinfeld}
\setcounter{equation}{0}

\subsection{Associateurs de Drinfeld et associateurs di\'edraux}

Soit $\la \in \k$. Un associateur de Drinfeld $\varphi(A,B)$
tel que défini dans \cite{DRIN} est l'exponentielle d'une s\'erie de Lie en deux variables $A$ et $B$,
que l'on peut consid\'erer comme l'exponentielle d'une
s\'erie de Lie en $A$, $B$, $C$ avec relation $A+B+C=0$. Elle
est soumise aux relations
\begin{eqnarray}
\varphi(B,A) = \varphi(A,B)^{-1} \label{drin1} \\
e^{\la A} \varphi(B,A) e^{\la B} \varphi(C,B) e^{\la C} \varphi(A,C) =1
\label{drin2}
\end{eqnarray}
ainsi qu'\`a une troisi\`eme relation qui ne nous sera
pas utile, appel\'ee \'equation du
pentagone. A un tel associateur on peut associer des morphismes
$B \to \mathfrak{B}(\k)$ o\`u $B$ et $\mathfrak{B}(\k)$
sont associ\'es non plus \`a des groupes di\'edraux mais
aux groupes de Coxeter de type $A_n$ pour tout $n \geq 2$. Si
$n=2$, il n'est pas n\'ecessaire d'imposer \`a $\varphi$ de v\'erifier
l'\'equation du pentagone pour que l'application d\'efinie par Drinfeld
soit un morphisme. Comme $A_2 = I_2(3)$, nous pr\'ecisons ici
les rapports entre ces associateurs de Drinfeld et le cas $m=3$.

Soit $m = 3$, $\varphi$ un associateur de Drinfeld, $\bar{\varphi} = 
\varphi(t_0,t_1) \in \exp \mathfrak{g}'(\k)$,
$\Phi = (\bar{\varphi})^{-1} = \varphi(t_1,t_0)$. On a $s \om
\bullet t_0 = t_1$, $s \om
\bullet t_1 = t_0$, $s \om
\bullet t_2 = t_2$. On d\'eduit alors de (\ref{drin1}) que $\Phi$ v\'erifie
(\ref{impeq2}), et de (\ref{drin2}) que $\Phi$ v\'erifie (\ref{impeq3}).
L'ensemble des associateurs de Drinfeld associ\'es \`a $\la \in \k$
s'injecte donc dans $\Ass^{(3)}_{\la}(\k)$.

\subsection{Le groupe de Grothendieck-Teichm\"uller}

Drinfeld introduit d'autre part un groupe $GT(\k)$, dit
de Grothendieck-Teichm\"uller, form\'e des couples de la forme
$(\la, g(X,Y))$, o\`u $\la \in \k$ et $g(X,Y)$, \'el\'ement
de la compl\'etion pro-$\k$-unipotente du groupe libre sur $X$ et $Y$,
satisfait \`a
\begin{eqnarray}
g(X,Y) = g(Y,X)^{-1} \label{drin3} \\
g(X_3,X_1) X_3^{\frac{\la -1}{2}} g(X_2,X_3) X_2^{\frac{\la -1}{2}}
g(X_1,X_2) X_1^{\frac{\la -1}{2}} =1 \label{drin4} 
\end{eqnarray}
pour $X_1 X_2 X_3 = 1$ ainsi qu'\`a une autre \'equation, encore
dite du pentagone.
Le groupe $P$ est engendr\'e par $p_1 = \sigma^2$, $p_2 = \tau^2$,
$p_3 = \sigma^{-1} \tau^2 \sigma = \tau \sigma^2 \tau^{-1}$.
soumis \`a la seule relation $p_3 p_2 p_1 = Z$ avec
$Z = (\sigma \tau)^3$ central. On pose $Y_i = p_i Z^{\frac{-1}{3}}
\in P(\k)$, $f = g( \tau^2 Z^{\frac{-1}{3}}, \sigma^2 Z^{\frac{-1}{3}})
= g(Y_2,Y_1) \in P(\k)$, et $O = \sigma \tau \sigma = \tau \sigma \tau$.
Comme $O \bullet \sigma = \tau$, $O \bullet \tau = \sigma$ et
$O \bullet Z = Z$, on a $O \bullet f = f$ d'apr\`es (\ref{drin3}). On veut
montrer (\ref{gexp2}), c'est-\`a-dire $( (\tau \sigma)^{-1} \bullet J)
\sss =\ooo f$ avec $\ooo = Z^{\frac{\la -1}{2}}$,
$\sss = (p_1)^{\frac{\la -1}{2}}$, $\ttt = (p_2)^{\frac{\la -1}{2}}$
et $J = \sss f^{-1} \ttt (\tau \bullet f)$. C'est \'equivalent
\`a 
$$
(\sigma^{-1} \tau^{-1} \bullet ( \sss f^{-1}))(\sigma^{-1}
\bullet \ttt) (\sigma^{-1} \bullet f) \sss f^{-1} = Z^{\frac{\la -1}{2}}.
$$
On a $\sigma^{-1} \bullet \ttt = \sigma^{-1} O \bullet \sss
= \tau \sigma \bullet \sss$, $\sigma^{-1} \bullet f = \sigma^{-1} O
\bullet f^{-1} = \tau \sigma \bullet f^{-1}$,
donc il s'agit de montrer
$$
(\sigma^{-1} \tau^{-1} \bullet F)(\tau \sigma \bullet F) F = Z ^{\frac{\la -1}{2}}
$$
avec $F = \sss f^{-1} = (Z^{\frac{1}{3}})^{\frac{\la -1}{2}}
Y_1^{\frac{\la -1}{2}}g(Y_2,Y_1)$. Or on a
$$
\left\lbrace
\begin{array}{rclcrcl}
\tau \sigma \bullet \sigma^2 & = & \tau \sigma^2 \tau^{-1} & \mbox{ d'o\`u }
& \tau \sigma \bullet Y_1 & = & Y_3 \\
\tau \sigma \bullet \tau^2 & = & \sigma^2  & \mbox{ d'o\`u }
& \tau \sigma \bullet Y_2 & = & Y_1 \\
\sigma^{-1} \tau^{-1} \bullet \sigma^2 & = & \tau^2 & \mbox{ d'o\`u }
& \sigma^{-1} \tau^{-1} \bullet Y_1 & = & Y_2 \\
\sigma^{-1} \tau^{-1} \bullet \tau^2 & = & \sigma^{-1} \tau^2 \sigma & \mbox{ d'o\`u }
& \sigma^{-1} \tau^{-1} \bullet Y_2 & = & Y_3 \\
\end{array}
\right.
$$
Ainsi il s'agit de v\'erifier
$$
Y_2^{\frac{\la -1}{2}}g(Y_3,Y_2) Y_3^{\frac{\la -1}{2}}g(Y_1,Y_3) 
Y_1^{\frac{\la -1}{2}}g(Y_2,Y_1) = 1
$$
si $Y_3 Y_2 Y_1 = 1$, c'est-\`a-dire  
$$
Y_2^{\frac{\la -1}{2}}g(Y_1,Y_2) Y_1^{\frac{\la -1}{2}}g(Y_3,Y_1) 
Y_3^{\frac{\la -1}{2}}g(Y_2,Y_3) = 1
$$
si $Y_1 Y_2 Y_3 = 1$, ce qui est bien \'equivalent \`a (\ref{drin4}). Ainsi
le groupe de Grothendieck-Teichm\"uller s'injecte dans $G^{(3)}(\k)$.

\subsection{L'associateur transcendant et $\varphi_{KZ}$}

On \'etablit ici le lien entre l'associateur $\varphi_{KZ}$
que Drinfeld construit \`a partir de l'\'equation de Knizhnik-Zamolodchikov
et l'associateur transcendant $\Phi_0$ pour $m=3$.

On note $D_1,D_2,D_0$ les droites de $E = \R^3$ euclidien d\'efinies
par les \'equations $D_1 : x_2 = x_3$, $D_0 : x_1 = x_2$, $D_2 :
x_1 = x_3$, $V = \{ (x_1,x_2,x_3) \in E \ \mid \ x_1 + x_2 + x_3 = 0 \}$
et on pose $D'_i = D_i \cap V$. Une base de $V^*$ est form\'ee des
fonctions $x = x_1 - x_2$, $y = x_1 - x_3$. On munit $V$ de la base
duale.

Soit $\Omega$ la 1-forme correspondant \`a la connection KZ
sur $E \setminus D_0 \cup D_1 \cup D_2$,
$$
\Omega = t_0 \dd \log (x_1 - x_2) + t_1 \dd \log (x_2 - x_3)
+ t_2 \dd \log (x_1 - x_3)
$$
et $V' = V \setminus D'_0 \cup D'_1 \cup D'_2$. On a $t_0 + t_1 + t_2 =0$,
et
$$
\Omega_{|V'} = t_0 \frac{\dd x}{x} + t_1 \frac{ \dd y - \dd x}{y-x}
+ t_2 \frac{\dd y}{y}.
$$
Soit $\phi : V \to V$ donn\'e par la matrice $M = \left(
\begin{array}{cc} \frac{-1}{2} & 1 \\ \frac{\sqrt{3}}{2} & 0 \end{array}
\right)$. Les fonctions et formes diff\'erentielles sont transport\'ees
par $\phi_* : V^* \to V^*$ de matrice $ ( ^t M)^{-1}$. On a
$$
\phi_*(x) = \frac{2}{\sqrt{3}} y , \phi_*(y) =
\frac{1}{\sqrt{3}} (\sqrt{3} x + y), \phi_*(x-y)
= \frac{1}{\sqrt{3}} y -x.
$$
On pose $\theta = \pi/3$, $\theta_r = r \theta$ d'o\`u $\theta_0 = 0$,
$\theta_1 = \pi/3$, $\theta_2 = 2 \pi/3$,
$v_r = (\cos \theta_r, \sin \theta_r)$ et $\alpha_r(v) = \det(v_r,v)$. On
a
$$
\alpha_0 = y, \alpha_1 = -x \frac{\sqrt{3}}{2} + y \frac{1}{2},
\alpha_2 = -x \frac{\sqrt{3}}{2} -y \frac{1}{2}.
$$
\begin{figure}
\begin{center}
\includegraphics{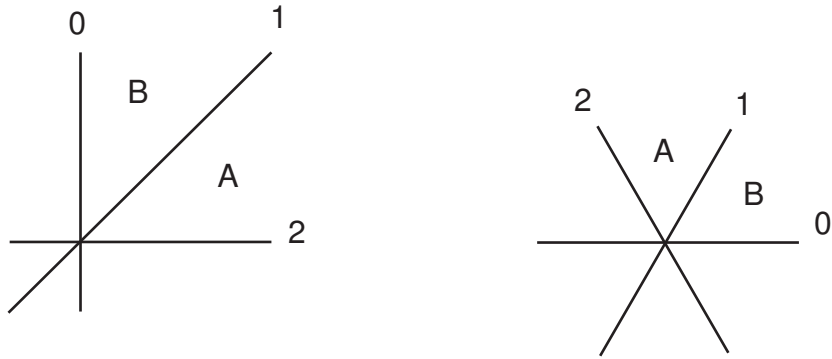}
\end{center}
\caption{Le changement de coordonn\'ees $\phi$} \label{figurephi}
\end{figure}
Notant $D''_i = \phi(D'_i) = \R v_i$, on a $\phi(V') = V''
= V \setminus D''_0 \cup D''_1 \cup D''_2$ (voir figure \ref{figurephi}) et
$$
\phi_*(\Omega_{|V'}) = t_0 \frac{\dd \alpha_0}{\alpha_0} + 
t_1 \frac{\dd \alpha_1}{\alpha_1} + 
t_2 \frac{\dd \alpha_2}{\alpha_2}
$$
On note $\eps : V \to V$ l'application $(s,w) \mapsto (s w , w)$. On a
$$
\begin{array}{lcl}
\eps^{-1} (D'_0) & = & \{ (s,w) \in \R^2 \ \mid \ w=0 \mbox{ ou } s=0 \} \\
\eps^{-1} (D'_1) & = & \{ (s,w) \in \R^2 \ \mid \ w=0 \mbox{ ou } s=1 \} \\
\eps^{-1} (D'_2) & = & \{ (s,w) \in \R^2 \ \mid \ w=0  \} 
\end{array}
$$
\begin{figure}
\begin{center}
\includegraphics{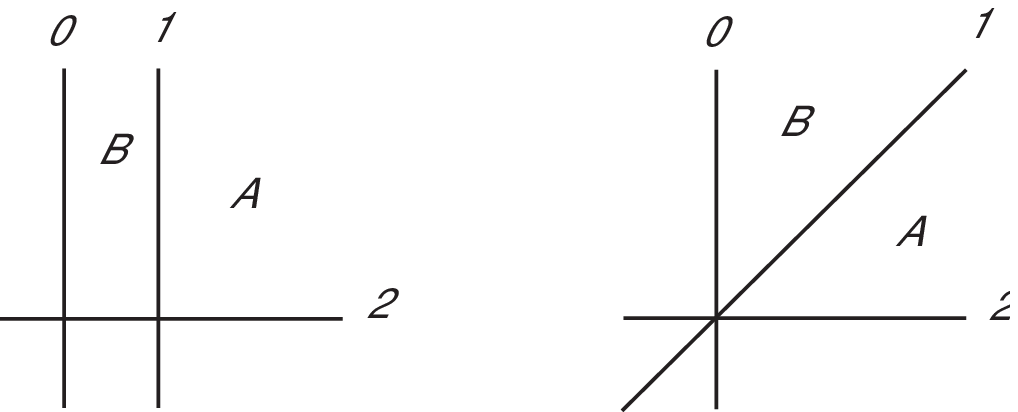}
\end{center}
\caption{L'\'eclatement $\eps$} \label{eclatement}
\end{figure}

On note $\Delta_0$, $\Delta_1$ et $\Delta_2$ les droites de $\R^2$
d'\'equations respectives $s=0$, $s=1$, $w=0$ et
$V^e = V \setminus \Delta_1 \cup \Delta_2 \cup \Delta_3$
(voir figure \ref{eclatement}). L'application
$\eps$ est un isomorphisme analytique $V^e \to V'$. De $x = s w$,
$y = w$ on d\'eduit
$$
\eps^*\left(\frac{\dd x}{x} \right) = \frac{\dd w}{w} + \frac{\dd s}{s}, 
\eps^*\left(\frac{\dd y}{y}\right) = \frac{\dd w}{w}, 
\eps^*\left(\frac{\dd x - \dd y}{x-y}\right) = \frac{\dd w}{w} + \frac{\dd s}{s-1}
$$
et
$$
\eps^* \left( \Omega_{|V'} \right) = \left( \frac{t_0}{s} + 
\frac{t_1}{s-1} \right) \dd s.
$$
L'\'equation diff\'erentielle $G'(s) = \left( \frac{t_0}{s}
+ \frac{t_1}{s-1} \right) G(s)$
sur $]0,1[$ admet deux solutions $G_{\pm}(s)$ d\'etermin\'ees
par $G_+(s) \sim (1-s)^{t_1}$ pour $s \to 1^-$, 
$G_-(s) \sim s^{t_0}$ pour $s \to 0^+$, et on a par d\'efinition
$\varphi_{KZ} = G_+^{-1} G_-$, c'est-\`a-dire $G_- = G_+
\varphi_{KZ}$. On en d\'eduit deux fonctions $F_{\pm} = G_{\pm}
(\frac{x}{y})$, solutions de $\dd F = \Omega F$ sur $\{ (x,y) \in V \ 
\mid \ 0<x<y \}$ uniquement d\'etermin\'ees par
$$
F_+ \sim \left( 1 - \frac{x}{y} \right)^{t_1} \mbox{ quand }
\frac{x}{y} \to 1^-, \ \  
F_- \sim \left( \frac{x}{y} \right)^{t_0} \mbox{ quand }
\frac{x}{y} \to 0^+.
$$
Or $\phi_*(\frac{x}{y}) = \frac{-\alpha_0}{\alpha_2}$, $\phi_*
(1 - \frac{x}{y}) = \frac{\alpha_1}{\alpha_2}$. On en d\'eduit
deux fonctions $\tilde{F}_{\pm}$ solutions de $\dd F = 
(\phi_* \Omega) F$ telles que
$$
\tilde{F}_+ \sim \left( \frac{\alpha_1}{\alpha_2} \right)^{t_1} \mbox{ quand }
\frac{-\alpha_0}{\alpha_2} \to 1^-, \ \  
\tilde{F}_- \sim \left( \frac{-\alpha_0}{\alpha_2} \right)^{t_0} \mbox{ quand }
\frac{-\alpha_0}{\alpha_2} \to 0^+.
$$
Passant en coordonn\'ees polaires, on a 
$$
\tilde{F}_+ \sim \left( \frac{\theta-u}{\sin(\theta)} \right)^{t_1} \mbox{ quand }
u \to \theta^-, \ \  
\tilde{F}_- \sim \left( \frac{u}{\sin(2\theta)} \right)^{t_0} \mbox{ quand }
u \to 0^+.
$$
Posons $\beta = 1/\sin \theta = 1/ \sin 2 \theta$. On a donc
$\tilde{F}_+ \beta^{-t_1} \sin (\theta -u)^{t_1}$ quand
$u \to \theta^-$, $\tilde{F}_- \beta^{-t_0} \sim u^{t_0}$
quand $u \to 0^+$, c'est-\`a-dire $\tilde{F}_+ \beta^{-t_1} = F_{1,-}$
et $\tilde{F}_- \beta^{-t_0} = F_{0,+}$. Comme $F_{1,-} =
F_{0,+} \Phi_0$ et $\tilde{F}_+ \varphi_{KZ} = \tilde{F}_-$, on a
$$
\varphi_{KZ} = \beta^{-t_1} \Phi_0^{-1} \beta^{t_0}
$$

\section{Appendice 3 : Associateurs d'Enriquez et morphismes en type $B$}

Dans \cite{ENRIQUEZ}, B. Enriquez définit des analogues des
associateurs de Drinfeld, qui lui permettent d'obtenir des morphismes
du groupe d'Artin de type $B_n$ vers les éléments grouplike
de certaines algèbres de Hopf complètes. Dans cet appendice nous établissons le
lien entre ces morphismes et ceux qui nous intéressent ici.

On désignera par W(S) le groupe de Coxeter associé au diagramme
de Coxeter $S$,
$A(S)$ (resp. $P(S)$) le groupe d'Artin (resp. le groupe
de tresses pures généralisé) associé. D'autre part, on note $\mathcal{B}_n
= A(A_{n-1})$ le groupe de tresses habituel, et $\mathcal{P}_n
= P(A_{n-1})$ le groupe de tresses pures.

\subsection{Le groupe d'Artin de type $B$ et les groupes $G(N,1,n)$}

On note $\tau, \sigma_1, \dots, \sigma_{n-1}$ les générateurs
standard de $A(B_{n})$. On a $\tau \sigma_1 \tau \sigma_1
= \sigma_1 \tau \sigma_1 \tau$, $\tau$ commute à $\sigma_i$ si $i >1$,
et les $\sigma_1,\dots, \sigma_{n-1}$ vérifient les relations de
$\mathcal{B}_n$. Ce groupe apparaît comme groupe de tresses généralisé
associé non seulement à $B_n$, mais à une famille de groupes
de réflexions complexes.

On note $W_{N,n}$ le groupe de réflexions complexes, de
type $G(N,1,n)$ dans la classification de ST, défini comme 
l'ensemble des éléments de $GL_N(\C)$ notés
$[\underline{a},\sigma]$ avec $\underline{a} = (a_1,\dots,a_n)$,
$a_i \in \C$ et $a_i^N = 1$, $\sigma \in \mathfrak{S}_n$,
tels que 
$$
[\underline{a},\sigma](z_1,\dots,z_n)= (a_1z_{\sigma(1)},\dots,a_nz_{\sigma(n)})
$$
Pour $N>1$ l'action de ce groupe est irréductible, et on a $W_{1,n}=
\mathfrak{S}_n$, $W_{2,n} = W(B_n)$.
On note
$s = [(\zeta,1,\dots,1), 1]$ avec $\zeta = e^{2 \ii \pi/m}$.
Dans tous les cas on a une décomposition de $W_{N,n}$ en produit semi-direct
$(\Z/N \Z)^n \rtimes \mathfrak{S}_n$, $s$ correspondant
à $(1,0,\dots,0) \in (\Z/N \Z)^n$.

Pour les propriétés de ces groupes et de leurs groupes de tresses
généralisé, nous renvoyons à \cite{BMR}. Rappelons que $s$
et les transpositions vérifient, outre les relations du groupe
symétrique sur les transpositions $s_i = (i \ i+1)$, les relations
$s s_1 s s_1 = s_1 s s_1 s$, $s^N=1$,  et $s$ commute à $s_i$
pour $i > 1$. D'autre part le complément dans $\C^n$ des hyperplans
de réflexions est
$$
X_{N,n} = \{ (z_1,\dots,z_n) \ \mid \ z_j \neq 0, z_j \neq \zeta^a z_k \}
$$
dont le groupe fondamental est noté $P_{N,n}$. On vérifie facilement que l'algèbre d'holonomie $\mathfrak{g}_{N,n} = \mathfrak{g}_{X_{N,n}}$
de $X_{N,n}$ est exactement l'algèbre $t_{n+1,N}$ de \cite{ENRIQUEZ}.
Cette dernière est engendrée par des générateurs $t(a)^{i,j}$ et $t_0^{1,i}$
pour
$i,j \in [2,n]$, et on trouve parmi les relations de définition
$t(a)^{ij} = t(-a)^{ji}$. Plus précisément, notant
$v(a)^{i,j} = t(a)^{i+1,j+1}$, $u_i =
t_0^{1,i+1}$, les générateurs
$v(a)^{i,j}$ correspondent aux hyperplans $z_i = \zeta^a z_j$, 
les éléments $u_i$ aux hyperplans
$z_i = 0$. Il suffit alors de vérifier que les relations de définition
sont les mêmes pour en conclure $\mathfrak{g}_{N,n} \simeq
t_{n+1,N}$ comme algèbre de Lie graduée.

Pour tout $N > 1$, le quotient $X_{N,n}/W_{N,n}$ a $A(B_n)$
comme groupe fondamental. La preuve du theorème 7.1 et
le corollaire 7.5 de \cite{ENRIQUEZ} montre l'existence
de morphismes $A(B_n) \to (\U \g_{N,n}(\k))^{\times}
 \rtimes
W_{N,n}$ pour tout corps $\k$ de caractéristique 0, donnés par
$$\tau \mapsto s \exp(u_1), \ \ \ \sigma_i = \Psi_i s_i \exp(N v(0)^{i,i+1}/2) \Psi_i^{-1} $$
où $\Psi_i \in (\U \g_{N,n}(\k))^{\times}$. Dans les
termes de l'introduction, pour $N=2$ cela montre
\begin{theore} Si $W$ est de type $B_n$, pour tout corps
$\k$ de caractéristique 0 il existe un morphisme $\PPhi : \k B \to \mathfrak{B}
(\k)$ qui vérifie la condition fondamentale.
\end{theore}
Pour $N>2$, cela montre l'existence de morphismes de
ce type pour les groupes de réflexions complexes $W_{N,n}$ et leur groupe
de tresses généralisé au sens de \cite{BMR}. Cela permet également
d'espérer élargir ces propriétés aux groupes de réflexions complexes.

\subsection{Les groupes $K_{N,n}$ et $P_{N,n}$}

B. Enriquez obtient des isomorphismes de la complétion pro-$\k$-unipotente
d'un sous-groupe $K_{N,n}$ de $\mathcal{P}_n$ avec
le groupe des éléments grouplike de $\U \g_{N,n}(\k)$.
Nous montrons que
ce sous-groupe $K_{N,n}$ est naturellement isomorphe à $P_{N,n}$.

Le groupe $K_{N,n}$ est défini dans \cite{ENRIQUEZ} comme suit.
On note $x_{ij}$ les générateurs classiques de $\mathcal{P}_n$,
et $F'_n$ le sous-groupe engendré par les $x_{1,i}$, $i \in [2,n]$. Il
est classique que $F'_n$ est libre sur ces générateurs, et que l'on a
une décomposition
$\mathcal{P}_n \simeq \mathcal{P}'_{n} \ltimes
F'_{n-1}$ où $\mathcal{P}'_n \simeq \mathcal{P}_{n-1}$ laisse fixe le
premier brin.
Il est classique que
l'action par conjugaison de $A(A_{n-1}) = \mathcal{B}_n$  
sur les générateurs libres $x_{1i}$ de $F'_n$
est l'action d'Artin. En particulier, $\mathcal{P}'_n$ agit
trivialement sur l'abélianisé $\Z^{n-1}$ de $F'_n$, et on en déduit un
morphisme $\mathcal{P}_n \to \Z^{n-1}$ puis par réduction modulo $N$
un morphisme $\mathcal{P}_n \to (\Z/N \Z)^{n-1}$, dont
$K_{N,n}$ est par définition le noyau. Dans \cite{ENRIQUEZ}, $\mathcal{P}_n$
(et donc $K_{N,n}$) est considéré comme un sous-groupe
de $A(B_{n-1})$. Cette interprétation, visuellement évidente
en termes de tresses géométriques, peut se décrire algébriquement
comme suit.

On a un morphisme
surjectif $A(B_n) \to A(A_{n-1})$ donné par $\tau \mapsto 1$,
$\sigma_i \mapsto \sigma_{i}$. Ce morphisme est scindé, par
$\sigma_i \mapsto \sigma_{i}$. Son noyau $F_n$ est
un groupe libre de rang $n$ sur les générateurs 
$$
\gamma_i = \sigma_{i-1} \sigma_{i-2} \dots \sigma_1 \tau \sigma_1^{-1}
\dots \sigma_{i-2}^{-1} \sigma_{i-1}^{-1}
$$
De plus, l'action de $A(A_{n-1})$ sur lui
est l'action d'Artin (cf. \cite{PARIS} prop 2.1) :
$\sigma_i$ envoie $\gamma_i$ sur $\gamma_{i+1}$,
$\gamma_{i+1}$ sur $\gamma_{i+1}^{-1} \gamma_{i} \gamma_{i+1}$ et
$\gamma_j$ sur $\gamma_j$ si $j \not\in \{i,i+1\}$.
En particulier, le morphisme d'abélianisation
$F_n \to \Z^n$ est équivariant, par rapport à l'action de $A(A_{n-1})$
sur $F_n$, l'action par permutation de $\mathfrak{S}_n$
sur $\Z^n$, et le morphisme naturel $A(A_{n-1}) \to \mathfrak{S}_n$.
On en déduit un morphisme surjectif $A(B_n) \simeq
A(A_{n-1}) \ltimes F_n \to \mathfrak{S}_n \ltimes \Z^n$. Il envoie
$\tau$ sur $(1,0,\dots,0) \in \Z^n$ et $\sigma_i$ sur $(i-1,i)
\in \mathfrak{S}_n$ pour $i > 1$. En composant avec la
projection canonique $\Z \mapsto \Z/N \Z$, on en déduit
des morphismes $\pi_N : A(B_n) \to \mathfrak{S}_n \ltimes (\Z/N \Z)^n$.
On a $\pi_N(\tau) = ( \bar{1},0,\dots,0) \in (\Z/N \Z)^n$,
$\pi_N(\sigma_i) = (i-1,i) \in \mathfrak{S}_n$ pour $i >1$. Il s'agit donc
du morphisme déjà mentionné $A(B_n) \to W_{N,n}$, en particulier
pour $N=2$
de la projection naturelle du groupe de tresses généralisé $A(B_n)$
sur le groupe de Coxeter $W(B_n)$.

On a alors $A(B_n) \simeq \mathcal{B}_n \ltimes F_n$. Considérons
l'isomorphisme $F'_n \to F_n$ défini par $x_{1i} \mapsto \gamma_{i-1}$.
On déduit de la compatibilité des actions le diagramme commutatif suivant,
dans lequel les flèches horizontales sont des isomorphismes
et les flèches verticales des injections.

$$
\xymatrix{
 & F'_n \rtimes A(A_{n-1}) \ar[r] & F_n \rtimes A( A_{n-1})  \ar[r] &  A(B_n)\\
\mathcal{P}_{n+1} \ar[r] & F'_n \rtimes \mathcal{P}_n \ar[r] \ar[u] & 
F_n \rtimes \mathcal{P}_n \ar[u]\\
}
$$
L'inclusion $\mathcal{P}_{n+1} \hookrightarrow A(B_n)$ est celle donnée
par la composition des flèches du diagramme.
Comme la restriction de $\pi_N$ à $F_n \rtimes \mathcal{P}_n$
est simplement $(F_n \to (\Z/N \Z)^n) \times 1$,
on en déduit bien que son noyau $P_{N,n}$ s'identifie à $K_{N,n+1}$.

\end{document}